\documentclass[11pt]{article}
\usepackage{latexsym}
\usepackage{amssymb}
\newtheorem{theorem}{Theorem}[subsection]

\newtheorem{corollary}[theorem]{Corollary}
\newtheorem{definition}[theorem]{Definition}
\newtheorem{example}[theorem]{Example}
\newtheorem{lemma}[theorem]{Lemma}

\newtheorem{proposition}[theorem]{Proposition}
\newtheorem{remark}[theorem]{Remark}

\newenvironment{proof}[1][Proof]{\textit{#1.} }{\ \hfill $\Box$}
\font\gordas = msbm10 at 12pt
\def\bbb#1{\hbox {{\gordas #1}}}
\def\erre{{\bbb R}}
\def\zet{{\bbb Z}}
\def\ene{{\bbb N}}
\def\UNO{1\mkern-7mu1}
\def\ee{{\bbb E}}
\def\pee{{\bbb P}}
\def\ene{{\bbb N}}

\begin{document}
\begin{center}
{\bf DEGREES OF TRANSIENCE AND RECURRENCE \\ AND HIERARCHICAL RANDOM WALKS}\\[.5cm]
\end{center}
\noindent
D.~A.~DAWSON$^1$\\
 {\it Carleton University}\\
{\it Ottawa, Canada K1S 5B6  (e.mail: ddawson@math.carleton.ca)}\\[.5cm]
L.~G.~GOROSTIZA$^2$ \\
{\it Centro de Investigaci\'on y de Estudios Avanzados}\\
{\it 07000 M\'exico D.F., Mexico (e.mail: lgorosti@math.cinvestav.mx)}\\[.5cm]
A.~WAKOLBINGER$^3$\\
 {\it Goethe-Universit\"at}\\
{\it Frankfurt am Main, Germany (e.mail: wakolbin@math.uni-frankfurt.de)}
\vglue .5cm
\centerline{(December 20, 2003)}
\vglue1cm
\noindent
\begin{abstract}  The notion of degree  and related notions concerning recurrence and transience for a class of L\'evy processes on
metric Abelian groups are studied. The case of random walks on a hierarchical group is examined with emphasis on the role of the ultrametric structure of the group and on analogies and differences with Euclidean random walks. Applications to separation of time scales and  occupation times of multilevel branching systems are discussed.
\end{abstract}
{\bf Mathematics Subject Classifications (2000):} 60J15, 60J30, 60B15,  60F05, 60J80.
\vglue1cm
\noindent
{\bf Key words:} degree, degree of transience, degree of recurrence,
 $k$-strong transience,
hierarchical group,
hierarchical random walk,  ultrametric space, separation of time scales,
 multilevel branching system, occupation time.

\footnote{\kern-.6cm $^1$ Research supported by NSERC (Canada) and a Max Planck Award for International Cooperation.\\
$^2$ Research supported by CONACYT grant 37130-E (Mexico).\\
$^3$ Research supported by  DFG (SPP 1033) (Germany).}

\tableofcontents
\vglue .5cm
\noindent
\section{Introduction}\label{sec1}
For a  L\'evy process $X$ on a metric Abelian group, we introduce its {\it degree} $\gamma$  as the supremum over all $\zeta > -1$ for which the operator power $G^{\zeta +1}$ of the Green operator of $X$ is finite (in a sense made precise in subsection \ref{GOP}). If $\gamma $ is positive, we call it the {\em degree of transience}, if it is negative, we call $-\gamma$ the {\em degree of recurrence} of $X$.
This extends notions defined in \cite{[9]}. For a transient process $X$, the degree of transience can be characterized as the order up to which the moments of last exit times $L_B$ (from a ball $B$ with positive radius) of $X$ are finite. For a large class of recurrent random walks on countable state spaces (at least for those whose transition probabilities $p_t$ have a power asymptotics in $t$) the degree of recurrence equals the order up to which the moments of first return times to the origin are finite.   

We say that $X$ has degree $\gamma^+$ (or alternatively, $\gamma^-$) if $X$ has degree $\gamma$ and 
$G^{\gamma +1}$ is finite (or infinite). For example, $d$-dimensional Brownian motion has degree $\gamma^-$ with $\gamma = d/2-1$.

A definition   of {\it $k$-strong transience}  for each integer $k\geq 1$
 was given in \cite{[9]}, which can be rephrased as follows: $X$ is $k$-strongly transient if the degree is bigger than $k$ or equal to $k^+$.  (The case $k=1$  corresponds to the usual strong transience.)   If the degree is either $k^-$ or $k^+$ then we say that $X$ is at the {\em border} of  $k$-strong transience (if $k \in \mathbb N$), or the border between transience and recurrence (if $k=0$). A process with degree $k^-$ is not $k$-strongly transient and for such a process the $(k+1)$-st operator power $G_t^{k+1}$ of the {\em incomplete potential operator} $G_t=\int_0^tT_s\,ds$ (where $T_t$ is the semigroup of the motion) typically has a subalgebraic growth as $t\to \infty$. More generally, we 
will define  operators $G_t^{({\zeta +1})}$, $\zeta>-1$, which in a certain sense interpolate between the integer powers of $G_t$, and we will investigate the growth of $G_t^{(\gamma +1)}$ as $t\to \infty$ for various examples of processes with degree $\gamma^-$.

A versatile class of random walks which (i) covers the range $(-1,\infty)$ of degrees, (ii) contains a wealth of examples with degrees $\gamma^-$ and $\gamma^+$, and (iii) allows a thorough analysis of cases on the borders, are
random walks on hierarchical groups (called {\it  hierarchical random walks}). They have their origin in the ``light bulb'' random walk studied by Spitzer \cite{[52]}
(page 93), and a model introduced by  Sawyer and Felsenstein \cite{[49]} in the context of genetics. For more background and references on hierarchical random walk we refer to the survey article \cite{[11]}.

The state space of the hierarchical walks is  $\Omega_N$, the {\it hierarchical group of
order} $N$. This is a countable Abelian  group consisting  of sequences of numbers in $\{0,1,\ldots,N-1\}$ only a finite number of which are different from zero,  with the metric such that the distance between two sequences is the largest coordinate number for which the respective coordinates are different. 

The countable group $\Omega_N$  is useful for the study of the large scale properties of hierarchical random walks. For the small scale properties it is necessary to pass to  a continuum hierarchical group (consisting of semi-infinite sequences). L\'evy processes on such a group have been considered in
\cite{[4], [19], [20], [38]}.

Hierarchical groups are examples of ultrametric spaces, where the distance $d(x,y)$ satisfies the strong (or non-archimedean) triangle inequality
$
d(x,y)\leq \max \{d(x,z), d(z,y)\}.
$
Ultrametric spaces are qualitatively different from Euclidean spaces; e.g.
two
balls  are either disjoint or
one contains the other.  Consequently,  a random walk on $\Omega_N$
can leave a closed ball of radius $R$ only by making a single jump of size greater than $R$ and not by a sequence of small jumps; in this respect, the hierarchical random 
random walks behave differently from  Euclidean 
random walks. On the other hand, some important aspects of the long-time behaviour ofÊ random walks as well as some classes of interacting
random walks depend only on their degree and consequently will be the same for random walks on Euclidean lattices and hierarchical groups of the same degree.
It is therefore of interest to calculate the degree of a Êrandom walk on a hierarchical group with jump distribution in a parametric family in terms of the parameters, and also to investigate
finer properties (such as the growth of $G^{(\gamma+1)}_t$ as $t\to\infty$) for a family of processesÊ all with degree $\gamma^{-}$.

For a given $N \ge 2$, $\mu >0$, and a sequence $(c_j)$ of 
positive numbers,  we consider the random walk on $\Omega_N$ which jumps distance $j$
with probability proportional to $c_{j-1}/N^{(j-1)/\mu}$, $j=1,2,...$, choosing  the arrival site 
uniformly among all sites at this distance. We call this  the
$(\mu, (c_j), N)$-random walk. The special case $c_j = 1$   gives the
$(\mu, (1), N)$-random walks, some of whose features were already
studied in \cite{[9]}. It is known that for a $(\mu, (1), N)$-random walk
$$p_t(0,0) \sim t^{-\mu}h(t),$$
where $h$ is a function which is bounded away from $0$ and $\infty$ and is ``slowly oscillating'' (i.e. $h(\log t)$ is periodic in $t$).
Consequently it has degree $ (\mu-1)^-$, and
the growth of $G^{(\mu)}_t$  is logarithmic as $t\to \infty$. For a
transient $(\mu, (1), N)$-random walk and $\zeta < \mu-1$, we
will study the asymtptotics of the last exit time moments 
$\mathbb E L_{B_R}^\zeta$ as $R\to \infty$, where $B_R$ is a ball
of radius $R$ containing the starting point.

Let us now turn to the more general $(\mu, (c_j), N)$-random
walks. We will show that under a mild condition on $(c_j)$, finiteness of $G^\mu$ is
equivalent to convergence of the sum $\sum_j c_j^{-\mu}$. For $\mu
=1$ and $\mu=2$ this amounts to transience and strong
transience of the walk. These summability conditions also play a
major role in connection with hierarchical equilibria of one- and
two-level branching populations \cite{[10]}.

With a view towards so-called mean-field limit  (see \cite{[10], [11]}),
it is of interest to study the behaviour of $(\mu, (c_j),
N)$-random walks as $N\to \infty$. It turns out that for a wide
class of sequences $c_j$, the degree of these random walks
approaches $\mu-1$ as $N\to \infty$. Indeed, we will show that for $0 < \liminf c_{j+1}/c_j
\le \limsup c_{j+1}/c_j < \infty$ the degree of the $(\mu, (c_j),
N)$-random walk is $\mu-1+O(1/\log N)$ as $N\to \infty$. If $\lim  c_{j+1}/c_j=1$, then the degree of
the $(\mu, (c_j), N)$-random walk is $\mu-1$ for {\em
all} $N$, and it is $(\mu-1)^+$ iff the $c_j^{-\mu}$ are summable.
For nondecreasing $c_j$ such that $\sum c_j^{-\mu}$ diverges,
the degree of the $(\mu, (c_j), N)$-random walk is
$(\mu-1)^-$, and $G^{(\mu)}_t$ grows like const $\sum_{j=0}^{\mu \log t/\log
N} c_j^{-\mu}$ as $t \to \infty$. In particular, the
$(\mu, ((j+1)^\beta), N)$-random walk (with $0<\beta$) has degree
$\mu-1$, and it has degree $(\mu-1)^-$ iff $\beta \le
\mu^{-1}$. In this case, $G^{(\mu)}_t$ grows like const $\log \log t$
for $\beta = \mu^{-1}$, and like const $(\log t)^{1-\beta \mu}$
for $0 < \beta < \mu^{-1}$. Proposition \ref{prop3.3.1} gives exact asymptotics for the growth of incomplete potential operator powers in some critical cases.

In subsection \ref{dmc} we investigate
the behaviour of the Markov chain given by the distance of a hierarchical walk to a
fixed point in $\Omega_N$, called {\it distance Markov chain}, which is the analogue of a (Euclidean) Bessel process. We will see that the distance Markov chains of hierarchical and Euclidean random walks behave differently, even if the underlying random walks have the same degree of transience/recurrence. For discrete time hierarchical random walks, we will study the distribution of the maximum of the distance Markov chain between times $0$ and $n$, and its asymptotics as $n \to \infty$, and we will see that for $\mu \ge 1$, $N^{j/\mu}$ is (asymptotically as $N \to \infty$) the right time scale for observing the exit behaviour of a $(\mu, (\eta^j), N)$- random walk from a closed ball of radius $j$. In \cite{[10]}
we study 2-level branching particle systems with a strongly transient migration on $\Omega_N$ which approaches the border of strong transience as $N\rightarrow \infty$, and which  leads to a separation of time scales and to a cascade of quasiequilibria associated to a sequence of nested balls of increasing radii,
in the  $N\rightarrow\infty$ limit. The results in subsection \ref{dmc} describe the appropriate time scale on each ball according to its radius and explain why, asymptotically as $N\rightarrow\infty$, only the evolution of the underlying random walk on the ball and on the surrounding ball of the next radius are relevant. This is the key for the cascade of quasiequilibria obtained in \cite{[10]} (see Remark 3.5.11). 

Part of the results obtained in the paper were motivated by  questions that 
arise in connection with
 occupation time fluctuations  and  hierarchical equilibria of branching 
systems studied in  \cite{[9], [10]}. The equilibrium behaviour and occupation time fluctuations of branching random walks and multilevel branching random walks provide examples of phenomena in which an essential role is played by the degree of the random walk, and also by the fine behaviour of the $(k+1)$-st powers of the incomplete potential operator $G_t$ when the random walk has degree $k^-$ with $k \in \mathbb N$. We will briefly review in Section 4 how the growth of $G_t^{k+1}$ and $GG_t^k$ as $t \to \infty$ carries over  to the growth functions in the norming of the occupation time fluctuations of $k$-level branching systems ($k=0,1,2$).
For $\alpha$-stable processes and $c^j$-random walks having degree $k^-$, it is known from \cite{[9]} that this growth function  is $\sqrt{t\log t}$. 
A case of particular interest, which is not covered by the results of \cite{[9]}, is provided by the $j^\beta$-random walks investigated in subsection \ref{sec3.3}. Here, we encounter a whole family of processes, all with degree $k^-$, leading to the (very) slow growth functions  $\sqrt{t(\log)^\delta}$, $0< \delta \le 1$, and $\sqrt{t\log \log t}$. 

Some of the results obtained in the paper  have been stated without proof in the survey article \cite{[11]}.

The paper is organized as follows.
Section 2 deals with  degree  and related notions,
 Section 3 refers to hierarchical random walks,
and Section 4 is devoted to occupation time fluctuations of branching systems.

\section{Degrees of transience and recurrence}
\subsection{Green operator powers and the degree of a L\'evy process}\label{GOP}
\setcounter{equation}{0}
 We consider L\'evy  processes  $X\equiv \{X(t), t\geq 0\}$ with cadlag paths  on $S$, a Polish space with (additive) Abelian group structure. We call $0 \,(\in S)$  the {\em origin} of $S$. For countable $S$, $X$ is a random walk on $S$ in continuous time.

The following function spaces will be used:
\vglue .15cm
\noindent
${\cal C}_b(S)$: continuous functions with bounded support,\\
${\cal B}_b(S)$: bounded measurable functions with bounded support,\\
${\cal C}^+_b(S), {\cal B}^+_b (S)$: elements of the previous spaces with
non-negative values.

Let $\{T_t, t\geq 0\}$ denote the semigroup of $X$, i.e., $T_t\varphi (x)=\ee_x \varphi (X(t)),\varphi \in {\cal B}_b(S)$. Recall that the potential (or Green) operator
 of $X$ is  defined by
$$
G\varphi =\int^\infty_0T_t \varphi dt, \quad \varphi \in {\cal B}_b (S),
$$
and  the fractional
 powers of $G$ are given by
\begin{equation}
\label{eq:2.1}
G^\zeta\varphi ={1\over \Gamma(\zeta)}\int^\infty_0 t^{\zeta-1}T_t \varphi dt,
\qquad \zeta>0,
\end{equation}
provided that the integrals are well defined. Note that
$G^{\zeta_1+\zeta_2}=G^{\zeta_1}(G^{\zeta_2})$ and for $\zeta=k$ integer, (\ref{eq:2.1}) coincides with the $k$th (operator) power of $G$. Recall that the process  $X$ is said to be {\it recurrent} iff $G\varphi \equiv \infty$ for
$\varphi \in {\cal C}^+_b (S)$, $\varphi \neq 0$, and {\it transient}
iff $||G\varphi|| <\infty$ for $\varphi \in {\cal C}^+_b(S)$
($||\cdot ||$ denotes the supremum norm).

\begin{definition}\label{def2.5}\rm
The degree $\gamma$ of $X$ is defined as
\begin{equation}\label{defdegree}
\gamma=\sup\{\zeta >-1: G^{\zeta+1}\varphi<\infty\quad\hbox{\rm for all}\quad\varphi\in {\cal B}^+_b(S)\}.
\end{equation}
 If $\gamma>0$, we call $\gamma$ the {\it degree of transience} of $X$, and if $-1<\gamma<0$, we call $-\gamma$ the {\it degree of recurrence} of $X$. The case $\gamma=0$ is considered in Definition \ref{def2.7}.
\end{definition}
\begin{remark}\label{rem2.6}\rm
\rm
 (a) These definitions extend Definition 2.4.2 in \cite{[9]} (see Remark \ref{rem2.10} below).\\
\noindent
(b) The definition of degree is valid without the Abelian group assumption, but in all the cases we consider here the space is an Abelian group and the processes are symmetric.
\end{remark}

In the transient case  we will relate $G^\zeta$ for $\zeta > 1$ to moments of last exit times (see subsection \ref{ssdegtrans}), and in the recurrent case we will relate 
$G^\zeta$ for $\zeta < 1$ to the finiteness of certain moments of first return times at least in special cases (see subsection \ref{ssdegrec}). 

If the degree  $\gamma$ defined by (\ref{defdegree}) is finite, it may be that
$G^{\gamma +1}\varphi<\infty$ or $G^{\gamma +1}\varphi\equiv \infty$, $\varphi \neq 0$.
In order to distinguish between the two cases and  abbreviate statements
 we introduce the following terminology:
\begin{definition}\label{def2.7}\rm
For a process $X$ of finite degree $\gamma$, we say that it has {\it degree $\gamma^+$} if
\begin{equation}
G^{\gamma +1}\varphi <\infty, \quad \varphi \neq 0,
\end{equation}
and it has {\it degree $\gamma^-$} if
\begin{equation}
G^{\gamma +1}\varphi \equiv \infty, \quad \varphi \neq 0.
\end{equation}
A process having degree 
$0^-$ will be called {\it critically recurrent}. 
\end{definition}

The symmetric $\alpha$-stable L\'evy process  on $\erre^d$, $0< \alpha \le 2$ (called $\alpha$-stable process henceforth) has degree $\gamma^-$ with
\begin{equation}\label{eq:2.19}
\gamma ={d\over \alpha}-1.
\end{equation}
We will also consider continuous time random walks on $\mathbb Z^d$  for which the jump distribution is in the domain of attraction of a symmetric $\alpha$-stable law and is $1$-lattice (i.e. the lattice generated by all vectors $x-y$ such that the transition probabilities $p_1(0,x)$ and $p_1(0,y)$ are strictly positive coincides with $\mathbb Z^d$). These walks will be called $(\alpha, d)$-{\em random walks} for short, and they also have degree $\gamma^-$ with $\gamma$ given by (\ref{eq:2.19}).     Indeed, combining a multidimensional local limit theorem (\cite{Rv}, Theorem 6.1) with a moderate deviations argument for Poisson random variables it is easy to see that the transition probability $p_t$ of an $(\alpha, d)$-random walk satisfies
\begin{equation}\label{loclim}
p_t(0,0) \sim  \mbox {const } t^{-d/\alpha} \mbox{ as } t \to \infty.
\end{equation}

Note that within the class of symmetric $\alpha$-stable processes and of $(\alpha, d)$-random walks the degree $\gamma$ is restricted to $[-1/2,\infty)$. Obviously, these processes  are critically recurrent for $d=\alpha$.
 For Brownian motion ($\alpha = 2$) on $\erre^d$ and simple symmetric random walk on $\zet^d$ the degree is $d/2-1$.
By the scaling property of the $\alpha$-stable process,
 $\ee_0 L^\zeta_{B_R}=R^{\alpha \zeta}\,\ee_0 L^\zeta_{B_1}$ for all $R$, where $\ee_0 L^\zeta_{B_1} < \infty$ for $\zeta <\gamma$. This growth  in $R$ willl be compared later on with corresponding results for certain hierarchical random walks.

A simple asymmetric random walk on $\zet$ has degree $\infty$. Also, Brownian motion on an infinite-dimensional Hilbert space (with nuclear covariance) has degree $\infty$.

Concluding this subsection, we recall the notions of strong/weak transience  and $k$-strong/weak transience \cite{[9]} which are closely related to the notion of degree and which play a role e.g. in connection with multilevel branching particle systems (see Section 4).

\begin{definition}\label{defstrongtrans}\rm
For each integer $k\geq 1$, we say that $X$ is
$$
\hbox{\it k-strongly transient}\;\; \;\hbox{\rm iff} \;\;||G^{k+1}\varphi ||<\infty \;\;\;\hbox{\rm for}\;\;\;\varphi \in {\mathcal B}^+_b (S),
$$
and
$$
k\hbox{\it -weakly transient}\;\;\;\hbox{\rm iff} \;\;\;||G^k\varphi ||<\infty \;\;\;\hbox{\rm for}\;\;\;\varphi \in {\mathcal B}^+_b (S)\;\;
$$%
$$\hbox{\rm and}\;\;
G^{k+1}\varphi \equiv \infty\;\;\;\hbox{\rm for}\;\;\; \varphi \in {\cal C}^+_b (S),\;\;\; \varphi \neq 0.
$$
\end{definition}

The case $k=1$ corresponds to the usual strong and weak transience.
Definition \ref{defstrongtrans}  is
compatible with (and more streamlined than) the one in \cite{[9]}
(Definition 2.1.1). In \cite{[9]} we referred to $G^{k+1}\varphi\equiv\infty$ as ``level $k$ recurrence'' because it corresponds to recurrence of ``level $k$ clans'' in  branching systems.
Note that
$k$-strong transience implies $\gamma\geq k$, and  $k$-weak transience implies $\gamma\in[k-1,k]$. Conversely, 
$\gamma > k$ implies $k$-strong transience, and  $\gamma\in (k-1,k)$ implies $k$-weak transience. We shall see in examples that certain critical behaviours occur when $\gamma$ takes an integer value. The $\alpha$-stable process is $k$-strongly transient iff $\alpha < d/(k+1)$ and $k$-weakly transient iff $d/(k+1) \le \alpha < d/k$.

\subsection{Degree of transience and moments of last exit times}\label{ssdegtrans}
\setcounter{equation}{0}

In this subsection we give a connection between the operator powers $G^\zeta$, $\zeta \ge 1$, defined in (\ref{eq:2.1}) and moments of last exit times. Intuitively, this relates to the degree of transience as follows:  the higher the degree of transience, the quicker the process tends to leave a bounded set forever. 

For a non-empty  Borel set $A\subset S$, let $L_A$ denote the last exit time 
 of $X$ from $A$,
$$L_A=
\sup \{t\geq 0: X(t)\in A\}\quad{\rm (if}\,\,\{t\geq0:X(t)\in A\}\neq \phi).
$$

\begin{proposition}\label{prop2.2}
Assume $X$ is transient  and  for any
closed ball $K\subset S$,
\begin{equation}
\label{eq:2.3}
\sup_{x\in K} G\UNO_K (x)<\infty,
\end{equation}
and for any closed ball $C\subset K^\circ$ (interior of $K$),
\begin{equation}
\label{eq:2.4}
\inf_{x\in C}G \UNO_K (x)>0.
\end{equation}
Then there exist positive constants $a_1$ and $a_2$ such that for all $\zeta >0$ and $x\in S$, 
\begin{equation}\label{eq:2.5}
a_1 G^{\zeta+1}\UNO_C(x)
\;\;\leq\;\; \ee_x L^\zeta_C\;\; \leq\;\;
a_2
G^{\zeta+1}\UNO_K(x).
\end{equation}
\end{proposition}

The  proof  is borrowed from  \cite{[47], [48]}.
 Those papers deal only with processes on $\erre^d$ 
 but the argument is general.\\
\begin{proof} Let $F_A=\inf \{t>0:X(t)\in A\}$ (the hitting time of $A\subset S$). By the Markov property of $X$  we have
\begin{eqnarray}
\label{eq:2.6}
G\UNO_K(x)&\geq & \ee_x \biggl(\UNO_{[F_C<\infty]}\ee_{X(F_C)}\int^\infty_0 \UNO_K (X(t))dt\biggr)\nonumber\\
&\geq & \inf_{y\in C}G \UNO_K(y)\pee_x (F_C <\infty).
\end{eqnarray}
By the Markov property and transience,
\begin{eqnarray}
\label{eq:2.7}
G\UNO_K(x)&=&\ee_x\biggl(\UNO_{[F_K<\infty]}\ee_{X(F_K)}\int^\infty_0 \UNO_K (X(t))dt\biggr)\nonumber\\
&\leq & \sup_{y\in K}G\UNO_K (y)\pee_x(F_K <\infty).
\end{eqnarray}
It follows from conditions (\ref{eq:2.3}) and (\ref{eq:2.4}), and from (\ref{eq:2.6}) and (\ref{eq:2.7}) that there exist positive constants $b_1$ and $b_2$ such that
\begin{equation}
\label{eq:2.8}
b_1 \pee_x (F_C <\infty)\;\;\leq \;\; G\UNO_K (x)\;\; \leq \;\;b_2 \pee_x (F_K<\infty)
\end{equation}
for all $x$.

Again by the Markov property,
\begin{equation}
\label{eq:2.9}
\ee_x L^\zeta_C =\int^\infty _0 \pee_x (L_C>t)\zeta t^{\zeta-1}dt
= \int^\infty_0 \ee_x \pee_{X(t)}(F_C<\infty)\zeta t^{\zeta -1}dt,
\end{equation}
since $L_C>t$ iff $F_C\circ \theta_t <\infty$, where $\theta _t$ is the shift of paths $\omega :(\theta_t \omega)(s)=\omega (t+s)$. Hence, by (\ref{eq:2.8}) and (\ref{eq:2.9}) there exist positive constants $b_3$ and $b_4$ such that
\begin{equation}
\label{eq:2.10}
b_3 \int^\infty_0 \ee_x G\UNO_C (X(t))\zeta t^{\zeta-1}dt\;\;\leq \;\;\ee_x L^\zeta_C\;\;\leq \;\; b_4 \int^\infty_0 \ee_x G\UNO_K(X(t))\zeta t^{\zeta -1}dt
\end{equation}
for all $x$.

Finally, for any closed ball $K$,
\begin{eqnarray}
\label{eq:2.11}
\lefteqn{
\int^\infty_0 \ee_x G\UNO_K (X(t))\zeta t^{\zeta-1}dt=\int^\infty_0 \ee_x
\int^\infty_0 \UNO_K (X(t+s))ds\zeta t^{\zeta -1}dt}\nonumber\\
&=&\int^\infty_0 \ee_x \int^\infty_t \UNO_K (X(s))ds \zeta t^{\zeta -1}dt=\ee_x \int^\infty_0 \UNO_K (X(s))\int^s_0 \zeta t^{\zeta -1}dtds\nonumber\\
&=& \int^\infty_0 s^\zeta T_s \UNO_K (x)ds,
\end{eqnarray}
and (\ref{eq:2.5}) follows from (\ref{eq:2.10}),  (\ref{eq:2.11}) and (\ref{eq:2.1}). 
\end{proof}
\smallskip

The following corollary is immediate.
\begin{corollary}
The degree of transience $\gamma\,\,(\geq0)$ is also given by
\begin{equation}\label{levelsets}
\gamma=\sup\{\zeta\geq 0:\ee L^\zeta_{B_R}<\infty\quad\hbox{\rm for all}\quad R>0\},
\end{equation}
where $B_R$ is a centered open ball of radius $R$. For irreducible transient  random walks on a countable Abelian group,
\begin{equation}\label{gammacount}
\gamma=\sup\{\zeta\geq 0,\, \ee L^\zeta<\infty\},
\end{equation}
where $L$ is the last exit time from $0$.
\end{corollary}

\begin{remark}\rm
For transient L\'evy processes on $\erre^d$, a set like on the r.h.s. of (\ref{levelsets})  is considered by Sato and Watanabe \cite{[47], [48]}.
\end{remark}

Conditions (\ref{eq:2.3}) and (\ref{eq:2.4}) hold in  all the examples considered in this paper.
\subsection{Degree of recurrence and moments of first return times}
\label{ssdegrec}
\setcounter{equation}{0}

In this subsection we only consider the case of countable $S$. We denote the transition probability of $X$ by $p_t(x,y)$. As before, we assume  that the walk $X$ is irreducible and (unless stated otherwise) starts in the origin.   
\begin{definition}\label{defT} \rm
Consider the {\em holding time}
$$H= \inf \{t>0: X_t \neq 0\}$$
and the {\em first return time} to the origin
$$T = \inf \{t>H: X_t = 0\}.$$
\end{definition}

For transient $X$ the last exit time $L$ from the origin is the sum of a geometric number of i.i.d. copies of $T$ {\em conditioned to be finite}, plus and independent copy of $H$.  Hence, in this case we have for all $\zeta >0$, 
$$\mathbb EL^\zeta < \infty \quad {\rm iff }\quad \mathbb E[T^{\zeta}\,| \,T<\infty] <
\infty.$$
Thus for transient $X$, the characterization (\ref{gammacount}) of the  degree of $X$ is equivalent to 
$$ \gamma = \sup\{\zeta\ge 0:\, \mathbb E[T^{\zeta}| T<\infty] <
\infty\}.$$

We now ask whether a similar characterization of the degree
in terms of moments of first return times also holds in the recurrent case . 
 
For the rest of this subsection we assume that $X$ is recurrent.
 Put $R= T-H$, and
$$\rho_t = \mathbb P[R>t],$$ that is, $1-\rho$ is the distribution function of the excursion time length $R$ of $X$ from the origin.
\medskip
\begin{lemma}\label{Spitzerlemma}
Assume rate $1$ holding times of $X$. Then for all $t>0$, 
\begin{equation}\label{Spitzercont}
\int_0^t p_s(0,0) \rho_{t-s}ds +p_t(0,0)= 1.
\end{equation}
\end{lemma}
\begin{proof}
We consider the process $Y_t := \UNO_{\{X_t \neq 0\}}$. The successive times $H_1 < H_2 < ...$ when $Y$ jumps from $0$ to $1$, together with the times $T_1 < T_2 < ... $ when $Y$ jumps back from $1$ to $0$, form an alternating renewal process, with the period in $0$ having distribution $\mathcal L(H)= \mbox{Exp}(1)$ and the period in $1$ having distribution
$ \mathcal L(R)$. Disintegrating the event $\{Y_t = 1\}$ with respect to the last jump of $Y$ from $0$ to $1$ before time $t$ we obtain
\begin{eqnarray*} \mathbb P[ Y_t = 1] &=& \int_0^t \mathbb P[H_i \in (s,s+ds), T_i > t \mbox { for some } i = 1,2,...]\\ &=&  \int_0^t \mathbb P[H_i \in (s,s+ds), T_i-H_i > t-s \mbox { for some } i = 1,2,...]\\
&=&  \int_0^t \mathbb P[H_i \in (s,s+ds) \mbox { for some } i = 1,2,...] \mathbb P[R > t-s] 
\\ &=&
\int_0^t p_s(0,0)\, \rho_{t-s}ds.
\end{eqnarray*}
The proof is complete since
$  \mathbb P[ Y_t = 0]= p_t(0,0)$.
\hfill \end{proof}
\smallskip
\begin{remark}\label{exactasympt}\rm

(a) Assume that for some $\mu > 0$ and a slowly varying function $\ell(t)$,
\begin{equation} \label{powerlaw}
p_t(0,0)\sim t^{-\mu}\ell(t) \mbox{ as } t \to \infty.
\end{equation}
Then the degree of the walk is $\gamma = \mu-1$. Indeed, (\ref{powerlaw})
implies that for each $\varepsilon > 0$, 
\begin{equation} \label{bounds}
p_t(0,0) \ge c_1 t^{-\mu-\varepsilon}  \quad {\rm and }\quad p_t(0,0) \le c_2 t^{-\mu+\varepsilon}
\end{equation}
for finite positive constants  $c_1, c_2$ depending on $\varepsilon$, and sufficiently large $t$.
Hence for all $\delta>0$, choosing $\varepsilon = \delta/2$ in (\ref{bounds}) we see that
$$\int_1^\infty t^{\mu-1+\delta} p_t(0,0) dt = \infty  \quad {\rm and }\quad  \int_1^\infty t^{\mu-1-\delta} p_t(0,0) dt < \infty.$$
The claim follows from (\ref{eq:2.1}) and (\ref{defdegree}).

(b) Assume that $p_t$ satisfies (\ref{powerlaw}) for $\mu \in (0,1)$ (as it is the case for $(\alpha, d)$-random walks with $d < \alpha$ and $\mu = d/\alpha$, see (\ref{loclim})). It follows from (\ref{Spitzercont})
that the Laplace transforms $\tilde p(\lambda)$ and $\tilde \rho (\lambda)$ of
$p_t(0,0)$ and $\rho_t$ are related by 
\begin{equation}\label{lap}
\tilde p(\lambda)\tilde \rho (\lambda) = \lambda^{-1} - \tilde p(\lambda),
\end{equation}
hence by a Tauberian theorem (\cite{[0]}, Theorem 1.7.6) one has
\begin{equation}
\tilde \rho(\lambda) \sim \lambda^{-\mu} \ell_1(1/\lambda) \quad \mbox{ as } \lambda \to 0
\end{equation}
for some slowly varying $\ell_1$.
Using another Tauberian theorem
 (\cite{[0]}, Theorem 1.7.2) one infers that
\begin{equation}\label{asymptrho}
\rho_t \sim t^{\mu-1} \ell_{\rho}(t) \mbox{ as } t \to \infty
\end{equation}
for some slowly varying function $\ell_{\rho}(t)$. 
Since
\begin{equation}\label{expR}
\mathbb E R^\zeta = \int_0^\infty \rho_t \zeta t^{\zeta-1} dt,
\end{equation}
we obtain from (\ref {asymptrho}), by a similar argument as in part (a), that
$$-\mu+1 = \sup\{\zeta \ge 0 : \,\mathbb E R^\zeta < \infty\}.$$ 
Since the first return time $T$ differs from $R$ only by the exponentially distributed holding time $H$,
and since the degree of the walk is $\gamma = \mu-1$  we have
\begin{equation}\label{chargammarec}
-\gamma = \sup\{\zeta\ge 0: \, \mathbb ET^{\zeta} <
\infty\}.
\end{equation}
\end{remark}

The next proposition shows that (\ref{chargammarec})  characterizes the degree of recurrence for all critically recurrent random walks satisfying the 
 the additional requirement
\begin{equation}\label{decayp}
p_t(0,0) = o(t^{-1+\varepsilon})\quad  \mbox{ as } t\to \infty \mbox{ for all } 0 < \varepsilon.
\end{equation}
Proposition \ref{prop3.2.4} and its corollary  show that an example of such a class of random walks are the $(1,(c_j),N)$-random walks (introduced in  Definition 3.1.4) where $(c_j)$ satisfies (3.2.16) and $\sum_j c_j^{-1} = \infty$.  
\begin{proposition}\label{nomoments}
For a recurrent random walk satisfying (\ref{decayp}), the return time $T$ has no moments of positive order.
\end{proposition}
\begin{proof}
We put $g_t= \int_0^t p_s(0,0)ds$.  For all $s, t > 0$ such that $p_r(0,0) < 1/2$ for all $r \ge s$, we have from (\ref{Spitzercont})
 $$1/2 \le 1-p_{s+t}(0,0) = \int_0^sp_r(0,0)\rho_{s+t-r} dr + \int_s^{s+t} p_r(0,0)\rho_{s+t-r} dr.$$
 Since $ \rho_t$ is decreasing, the first term on the r.h.s is bounded by
 $g_s \rho_t$, and the second one is bounded by $g_{s+t}-g_s$. Hence we obtain
 $$1/2 \le  g_s \rho_t + g_{t+s}-g_s.$$
 Using (\ref{decayp}), we have for each
 $0 < \varepsilon < 1$ and suitable constants $c_1, c_2 > 0$ depending on $\varepsilon$,
 \begin{eqnarray*}
 \rho_t \ge \frac{1/2-(g_{s+t}-g_s)}{g_s} 
 \ge \frac{1/2-c_1[(s+t)^\varepsilon-s^\varepsilon]}
 {c_1\,s^\varepsilon} 
 = c_2 s^{-\varepsilon}-((1+t/s)^\varepsilon - 1).
 \end{eqnarray*}
 Putting $s= t^{1/(1-2\varepsilon)}$ this turns
 into
 \begin{eqnarray*}
c_2\, t^{-\varepsilon/(1-2\varepsilon)}-
 \left(1+t^{-2\varepsilon/(1-2\varepsilon)}
 \right)^\varepsilon + 1 
 \sim c_2 \,t^{-\varepsilon/(1-2\varepsilon)}
 -\varepsilon t^{-2\varepsilon/(1-2\varepsilon)} \quad {\rm as } \quad t \to \infty.
 \end{eqnarray*}
 This shows that $\rho_t$ decays slower than $t^{-\delta}$ for
 any $\delta > 0$, and in view of (\ref{expR}) completes the proof.
 \end{proof}
 \smallskip
 
 The next proposition shows that for $\mu \in (0,1)$ a less restrictive condition than (\ref{powerlaw}) assures at least that the first return time has all moments of order less than $1-\mu$. This condition is fulfilled by the $(\mu, (c_j), N)$-random walks with $\mu \in (0,1)$ and  $(c_j)$ satisfying (3.2.16) (see Proposition \ref{prop3.2.4}).
 
 \begin{proposition} \label{prop2.2.5}
 For $\mu \in (0,1)$, assume 
\begin{equation}\label{lowerboundp}
p_t(0,0)^{-1} = o(t^{\mu+\varepsilon})\quad \mbox{ as } t \to \infty \mbox{ for all }\, 0< \varepsilon 
\end{equation}
(and consequently $\gamma \le \mu -1$).
Then the return time $T$ has all moments of order less than $1-\mu$.
 \end{proposition} 
\begin{proof}
Since
 $\rho_t$ is decreasing, we have
 from (\ref{Spitzercont})
\begin{equation}\label{boundrho}
1\ge 1-p_t(0,0) = \int_0^tp_s(0,0)\rho_{t-s} ds \ge \rho_t g_t.
\end{equation}
 From (\ref{lowerboundp}) we have that for each $\varepsilon > 0$ 
 there exists a constant $c>0$
 such that $p_t(0,0) \ge c t^{-\mu-\varepsilon}$, and consequently  
 $g_t\ge c_1t^{1-\mu-\varepsilon}$ for some $c_1 > 0$. Hence because of (\ref{boundrho}) we have for each $\varepsilon > 0$ and a suitable constant $c_\varepsilon $ 
 \begin{equation}\label{decayofrho}
 \rho_t 
 \le c_\varepsilon t^{\mu + \varepsilon -1}, \quad t> 0.
 \end{equation}
Consequently, for all $\delta \in (0, 1-\mu)$, putting $\varepsilon = \delta/2$ in (\ref{decayofrho}), we have from (\ref{expR})
$$\mathbb ER^{1-\mu-\delta} = \int_0^\infty  (1-\mu-\delta) t^{(1-\mu- \delta) -1}\rho_t dt \le 
{\rm const} \int_1^\infty  t^{-1-\delta/2}dt + {\rm const} < \infty .$$ 
Then it suffices to recall that $T = H+R$, where $H$ is exponentially distributed and therefore has moments of all orders.
\end{proof}
\smallskip
\begin{remark}\rm
Put $\gamma = \mu-1$. 

(a) For $\mu \in (0,1)$ the power asymptotics 
(\ref{powerlaw}) implies the equality (\ref{chargammarec}) (Remark \ref{exactasympt} (b)), which in this case characterizes the degree of recurrence in terms of moments of first return times (Remark \ref{exactasympt} (a)). 

(b) For $\mu =1$, the ``weak'' power asymptotics  (\ref{bounds}) (right part) still guarantees (\ref{chargammarec}), see Proposition \ref{nomoments}.  

 (c) For $\mu \in (0,1)$, the ``weak'' power asymptotics
(\ref{bounds}) (left part) implies that the return time $T$ has all moments of order less than $1-\mu = -\gamma$ (see Proposition \ref{prop2.2.5}). Hence in this case we have at least the bound 
\begin{equation}\label{onedirection}
-\gamma \le \sup\{\zeta\ge 0: \, \mathbb ET^{\zeta} <
\infty\}.
\end{equation}

(d) It would be interesting to know whether the characterization (\ref{chargammarec}) holds in general for recurrent random walks with degree $\gamma$.
\end{remark} 

\subsection{Incomplete potentials}\label{sec2.3}
\setcounter{equation}{0}
We now define the {\em incomplete potential operator} $G_t$ which together with its powers
 plays a basic role in occupation time results.
\begin{definition} \rm
For a process $X$ on $S$  we define the operator
\begin{equation}\label{Gt}
G_t\varphi=\int^t_0T_s\varphi ds,\quad \varphi\in {\cal B}_b(S),
\end{equation}
where $\{T_t\}$ is the semigroup of $X$. Moreover, we denote by 
$G^k_t,k=2,3,\ldots$ the  (operator) powers of $G_t$.
\end{definition}

When $G^k\varphi < \infty$,  $G^{k+1}_t\varphi\rightarrow\infty$ as $t\rightarrow\infty,\varphi\in
{\cal B}^+_b(S),\varphi\neq 0$,  the order of the growth of
$G^{k+1}_t\varphi$  
determines the appropriate normings for the occupation times in 
$k$-level branching populations. This is discussed in section 4.

For the  $\alpha$-stable  process on $\erre^d$ (having degree $\gamma = d/\alpha-1$)  and  integer $k \ge 0$ ,
\begin{equation}\label{eq:2.22}
\begin{array}{lll}
G^{k+1}_t\sim\kappa\log t&{\rm for}&\gamma=k \quad (\mbox{ equivalently } \alpha=\frac{d}{k+1}),\\ \nonumber
G^{k+1}_t\sim\kappa t^{k-\gamma}&{\rm for}&k-1<\gamma <k  \quad (\mbox{equivalently }
\frac{d}{k+1}<\alpha<\frac{d}{k})
\end{array}
\end{equation}
and
\begin{equation}\label{eq:2.23}
T_t \sim \kappa t^{-(\gamma +1)} \quad \mbox{as }  t\rightarrow\infty.
\end{equation}
In these formulas $\kappa$ stands for a positive constant which is different in each case, and formulas (\ref{eq:2.22}) and (\ref{eq:2.23}) are symbolic.
For example, the precise meaning of $G_t\sim \kappa t^{-\gamma}$ is
$\int_S\varphi G_t \psi d\rho\sim \kappa t^{-\gamma}H(\varphi, \psi)$, $\varphi, \psi \in {\cal C}_b (S)$, where $\rho$ is the Lebesgue measure on $\erre^d$ and $H(\varphi, \psi)$ is some positive-definite bilinear form
 \cite{[9]}. The ``critical'' cases $\gamma=k$ are associated with slowly varying
growth of $G^{k+1}_t$. 

We use the following notations:

 $a_t \asymp b_t$ as $t \to \infty$ \quad if $a_t/b_t$ and  $b_t/a_t$ remain bounded as $t \to \infty$, and

$a_t \propto b_t$ as $t \to \infty$ \quad if  $a_t/b_t$ and  $b_t/a_t$ are $o(t^\varepsilon)$ as $t \to \infty$ for all $\varepsilon > 0$.

The same notations will be used also for discrete indices $j=1,2,..$ in place of $t$.

\begin{remark}\label{rem2.10}\rm
(a) For  transient processes it is useful to consider the operator  $R_t$ defined by
$$
R_t =G-G_t, \quad t>0
$$
(see \cite{[42]}). 
It is easy to see that if $R_t\asymp t^{-\gamma}$ as $t\rightarrow\infty$ for some $\gamma >0$ (called transience of order $\gamma$ in \cite{[9]}), then the process is transient with  degree $\gamma^-$, and if $G_t\asymp t^{-\gamma}$ as $t\rightarrow\infty$ for some $\gamma\in (-1,0)$ (called recurrence of order $-\gamma$ in \cite{[9]}), then the process is recurrent with  degree $\gamma^-$. 

\noindent
(b) Recurrent processes such that $G_t=o(t^\varepsilon)$ as $t\rightarrow\infty$ for all $\varepsilon >0$ are critically recurrent. Indeed, for any $\zeta\in (-1,0)$ and  $0<\varepsilon <-\zeta$,
$$\int^\infty_1t^\zeta T_t\varphi dt\leq\sum^\infty_{k=0}(2^\zeta)^k\int^{2^{k+1}}_{2^k}T_t\varphi dt \leq {\rm const} \sum^\infty_{k=0}(2^\zeta)^k(2^{k+1})^\varepsilon <\infty .$$
An important case is $G_t\sim \mbox {\rm  const } \log t$ (which was called critical recurrence in \cite{[9]}).
\end{remark}

For processes with degree  $\gamma$, from the viewpoint of occupation times it is necessary to compute the growth of  $GG_t^{k-1}$ and of $G_t^{k}$  for an integer $k$ with
 $ \gamma +1 \le k < \gamma +2$ (cf. section \ref{sec4}).
However, also in the case of non-integer  $\zeta \ge \gamma+1$ it is interesting to study  the growth   of the operators
\begin{equation}\label{defGzetabrackets}
G^{(\zeta)}_t\varphi = \frac1{\Gamma(\zeta)}\int_0^ts^{\zeta-1} T_s\varphi \,ds, \quad t>0,\, \varphi \in \mathcal B_b(S)
\end{equation}
for $\varphi \ge 0$, $\varphi \neq 0$. Indeed,
the following lemma and its corollary show that for integer $k \ge \gamma +1$, and
a large class of walks, $G^{(k)}_t$ captures at least the growth of  $G^{k}_t$. Note that if $T_t\varphi \asymp t^{-(\gamma+1)}$, then 
\begin{equation} \label{asGtk}
G_t^{(k)} \asymp \left\{
\begin{array}{ll}
t^{k-(\gamma+1)}&{\rm if}\quad k>\gamma+1,\\
\log t&{\rm if}\quad  k=\gamma+1.
\end{array}\right.
\end{equation}
\medskip
\begin{lemma} For $k=1,2,...$ and $\varphi \in {\cal B}^+_b (S)$, $\varphi \neq 0$,\\
(a)
\begin{eqnarray} \label{sandwich}
 0\,\,\leq \,\,G^k_t\varphi - G^{(k)}_t\varphi\,  \leq \,\int^t_0ds_1\int^t_0ds_2\ldots \int^t_0ds_{k-1}
\int^{t+s_1+\ldots+s_{k-1}}_tT_s\varphi \, ds.
\end{eqnarray}
(b) Assume $T_t\varphi \asymp t^{-(\gamma+1)}$ for some $\gamma > -1$. Then for integer $k \ge \gamma +1$,
\begin{equation}\label{GtGt}
G_t^{(k)}\varphi \le G_t^{k}\varphi \le  G_t^{(k)}\varphi + {\rm const}\, t^{k-(\gamma+1)},
\end{equation}
and therefore
\begin{equation} \label{asGtkGtk}
G_t^{(k)}  \left\{
\begin{array}{ll}
\asymp G_t^{k}&{\rm if}\quad k>\gamma+1,\\
\sim G_t^{k}&{\rm if}\quad  k=\gamma+1.
\end{array}\right.
\end{equation}

\end{lemma}
\begin{proof}
For $k=1$, $G_t\varphi=G^{(1)}_t\varphi$, so there is nothing to prove.
For $k\geq 2$,
\begin{equation}\label{equalk}
G^{(k)}_t\varphi=\int^t_0T_sG^{(k-1)}_{t-s}\varphi \,ds,
\end{equation}
since from (\ref{defGzetabrackets}) the derivatives w.r. to $t$ of both sides of (\ref{equalk}) coincide by the semigroup property. 

Iterating  (\ref{equalk}), 
\begin{equation}\label{itint}
G^{(k)}_t\varphi=\int^t_0ds_1\int^{t-s_1}_0ds_2\ldots\int^{t-s_1-\ldots -s_{k-1}}_0ds_kT_{s_1+\ldots +s_k}\varphi.
\end{equation}
On the other hand, 
\begin{equation}\label{Gk}
G^k_t\varphi=\int^t_0ds_1\int^t_0ds_2\ldots\int^t_0ds_kT_{s_1+\ldots+s_k}\varphi,
\end{equation}
Subtracting  (\ref{itint}) from  (\ref{Gk}) we find
\begin{eqnarray}\label{sandwich2}
 0\,\,\leq \,\,G^k_t\varphi - G^{(k)}_t\varphi\, \leq\, \int^t_0ds_1\int^t_0ds_2\cdots\int^t_0ds_{k-1}\int^t_{t-s_1-\ldots -s_{k-1}}ds_kT_{s_1+\ldots +s_k}\varphi.
\end{eqnarray}
Substituting $s = s_1+...+s_k$ in the r.h.s. of (\ref{sandwich2}) we have 
(\ref{sandwich}).

Under the assumption of part (b), (\ref{GtGt}) is immediate from (\ref{sandwich}), and (\ref{asGtkGtk}) follows from (\ref{asGtk}) and
 (\ref{GtGt}).
\end{proof}
\begin{corollary}
For $\gamma > -1$, an integer $k \ge \gamma +1$ and $\varphi \ge 0$, $\varphi \neq 0$,  \\ \noindent
(a) if $T_t\varphi \asymp t^{-(\gamma +1)}$ then 
\begin{equation} \label{asGtt}
G_t^{k}\varphi \asymp \left\{
\begin{array}{ll}
t^{k-(\gamma+1)}&\mbox{ if}\quad k>\gamma+1,\\
\log t&\mbox{ if}\quad  k=\gamma+1,
\end{array}\right.
\end{equation}
\\ \noindent
(b) if  $T_t\varphi \,\propto\, t^{-(\gamma +1)}$,  then $G^k_t \varphi \, \propto\, t^{k-(\gamma +1)}.$ \end{corollary}

An example for (a) in the preceding corollary is given by the $(\gamma+1, (1), N)$-random walks (see (\ref{asp}) and Remark \ref{rem3.1.5}), and an example for (b) is provided by the $(\gamma+1, (c_j), N)$-random walks with $c_{j+1}/c_j \to 1$ (see Proposition \ref{prop3.2.4}).
\begin{definition}\rm
For discrete $S$ (as in the case of the hierarchical random walks studied in the following section) and $\zeta > 0$ we put 
\begin{equation}\label{defgtzeta}
g^{(\zeta)}_t= \frac 1{\Gamma(\zeta)}\int_0^ts^{\zeta-1} p_s(0,0) ds = G^{(\zeta)}_t\UNO_{\{0\}}(0), \quad t > 0, 
\end{equation}
\end{definition}
\begin{equation}\label{gzeta} 
g^{(\zeta)}= \frac 1{\Gamma(\zeta)}\int_0^\infty s^{\zeta-1} p_s(0,0) ds  = 
G^\zeta\UNO_{\{0\}}(0).
\end{equation}

For the $(\alpha, d)$-random walk (having degree $\gamma^-$ with $\gamma  =  d/\alpha-1$), (\ref{loclim}) implies
\begin{eqnarray}\label{gg} g^{(\gamma +1)}_t &\sim& \mbox{ const } \log t, \nonumber \\
g^{(\zeta)}_t &\sim& \mbox{ const }  t^{\zeta-(\gamma +1)} \quad \mbox{ for } \zeta > \gamma +1. \end{eqnarray}

\section {Random walks on the hierarchical group}\label{sec3}

\subsection {Hierarchical random walks}\label{sec3.1}
\setcounter{equation}{0}
\begin{definition}\label{def3.1.1}\rm
Let $N$ be an integer $\geq 2$.
The (countable) {\it hierarchical group of order $N$}
is defined by
$$
\Omega_N=\{x=(x_1,x_2,\ldots):x_i\in\zet_N,x_i\neq 0\,
{\rm only\, for\, finitely\, many}\, i\},$$
where $\zet_N=\{0,1,\ldots, N-1\}$ is the cyclic group of order $N$, with
 addition componentwise mod$(N)$. In other words, $\Omega_N$ is the direct sum of a countable number of copies of $\zet_N$. We endow
$\Omega_N$
 with the translation-invariant {\it hierarchical distance} $|\cdot |$ defined by
$$|x-y|=\left\{\begin{array}{lll}
0&{\rm if}&x=y,\\
\max\{i:x_i\neq y_i\} &{\rm if}&x\neq y.
\end{array}\right.$$%
\end{definition}
Note that
$d(x,y)=|x-y|$ is an ultrametric.
\begin{definition}\rm
\label{def3.1.2} ($r_j$-{\it random walk)}.
We consider
{\it hierarchical random walks} $\{\xi_n\}$ on  $\Omega_N$ defined by
$\xi_n=\sum\limits^n_{i=1}\eta_i, n=1,2,\ldots$, where
 $\eta_i, i=1,2,\ldots$ are   i.i.d. random variables in $\Omega_N$
 with distribution
\begin{equation}\label{3.1.1}
\mathbb P[\eta=y]={r_{|y|}\over N^{|y|-1}(N-1)},\quad y \in \Omega_N,
\quad  y\neq 0,\quad \mathbb P[\eta=0]=0,
\end{equation}
and $\{r_j,j=1,2,\ldots\}$ is a probability distribution on $\ene=\{1,2,\ldots\}$.
 Note that the random walk jumps from $x$ to $y$ such that $|x-y|=j\geq 1$ by first choosing distance $j$ with probability $r_j$ and then choosing $y$
uniformly on the sphere of radius $j$ with center at $x$
(since $N^{j-1}(N-1)$ is the number of points at distance $j$ from a given point). We assume that $r_j>0$ for all $j$, hence these random walks are irreducible.We call {\it $r_j$-random walk} the random walk defined by (3.1.1).  We will introduce descriptive names for special choices of $r_j$, and in some cases simplified names for easy identification; the name ${r}_j$-random walk always refers to the general case.
\end{definition}
\begin{remark}\rm
\label{rem3.1.3} (a)
The  $r_j$-random walks are the most general ``symmetric'' random walks on
$\Omega_N$ in the sense of  the uniform choice of a point at a given distance.

\noindent
(b) $\Omega_N$ can also be represented as the set of leaves
of a tree $T_N$. Each inner node of $T_N$ is at some level (or depth) $j \ge 1$, and the leaves are at level $0$. Each inner node at level $j$ has one neighbouring node at level $j+1$ (its {\em parent}) and $N$ neighbouring nodes at level $j-1$ (its {\em children}).  For a leaf $x$, let $a_j(x)$, $j=1,2,...$ denote its chain of ancestors. The $r_j$-random walk jumps from the leaf $x$ with probability $r_j$ to a leaf uniformly  chosen among all the leaves which descend from $a_j(x)$ but not from $a_{j-1}(x)$. 
\noindent
(c) The case $N=2$ corresponds to the ``light bulb'' random walk in  \cite{[52]}.
A criterion for transience/recurrence in this case was given in \cite{[7]}, and
extended in \cite{[22]} allowing  $N$ to
 depend on the index of each component.
Sawyer and Felsenstein \cite{[49]} used random walks  on $\Omega_N$  to
study genetic relatedness in a spatially structured population.
 It would be interesting to study hierarchical  random walks with random $N$
 (i.i.d. numbers of outgoing edges from each inner node), and non-symmetric hierarchical random walks.
\end{remark}

The $n$-step transition probability $p^{(n)}(x,y)$ of the $r_j$-random walk $\{\xi_n\}$, which can be obtained by Fourier methods \cite{[49], [23], [35]}, is given by
\begin{eqnarray}\label{3.1.2}
p^{(n)}(0,y)&=&(\delta_{0,|y|}-1){f^n_{|y|}\over N^{|y|}}
+(N-1)\sum^\infty_{k=|y|+1}{f^n_k\over N^k},
\quad n\geq 1,\quad y\in\Omega_N\backslash \{0\},\\ \nonumber
p^{(1)}(0,0)&=& 0,\kern9.7cm
\end{eqnarray}
where
\begin{equation}\label{3.1.3}
f_k=\sum ^{k-1}_{j=1}r_j-{r_k\over N-1}=1-r_k{N\over N-1}-\sum^\infty_{j=k+1}r_j,\quad k\geq 1.
\end{equation}

A continuous-time random walk $X=\{X(t),t\geq 0\}$ on $\Omega_N$   corresponding to  $\{\xi_n\}$, with unit  rate holding time, i.e., with transition probability
\begin{equation}\label{3.1.4}
p_t(0,y)=e^{-t}\sum^\infty_{n=0}{t^n\over n!}p^{(n)}(0,y),
\quad t\geq 0, \quad y\in\Omega_N,
\end{equation}
is given by
\begin{equation}\label{3.1.5}
p_t(0,y)=(\delta_{0,|y|}-1){e^{-h_{|y|}t}\over N^{|y|}}+(N-1)
\sum^\infty_{j=|y|+1}{e^{-h_jt}\over N^j},\quad t\geq 0,\quad y\in\Omega_N,
\end{equation}
where $h_j=1-f_j$, i.e.,
\begin{equation}\label{3.1.6}
h_j=r_j{N\over N-1}+\sum^\infty_{i=j+1}r_i,\,\, j=1,2,\ldots,
\end{equation}
\cite{[23], [35]}.  The $r_k$ are obtained from the $h_k$ by
\begin{equation}\label{3.1.7}
r_k={N-1\over N}h_k-{(N-1)^2\over N}N^k\sum^\infty_{j=k+1}{h_j\over N^j},
\quad k=1,2,\ldots .
\end{equation}
\begin{definition}\label{def3.1.4}\rm
 ($(\mu,(c_j),N)$-{\it random walk)}.
We  consider $r_j$-random walks   (3.1.1) with jump probabilities $r_j$
of the form
\begin{equation}\label{3.1.8}
r_j=D {c_{j-1}\over N^{{(j-1)}/\mu }},\quad j=1,2,\ldots ,
\end{equation}
where $\mu$ is a positive constant, $\{c_j,j=0,1,\ldots\}$ is a sequence of positive numbers and $D$ is a normalizing constant. This random walk as well as its continuous time version with unit rate holding times
will be called {\it $(\mu, (c_j), N)$-random walk,} emphasizing the three elements that define the jump probabilities. (It will be clear in each case whether the time is continuous or discrete.)
\end{definition}
\begin{remark} \label{rem3.1.5}\rm
For fixed $N$ and $\mu \neq 1$,  a  $(\mu,(\eta_j),N)$-random walk is the same  as a $(1,(c_j),N)$-random walk with $c_j=\eta_jN^{j(\mu-1)/\mu}$.
This transformation is useful because we also are interested in the behaviour of $(\mu,(c_j),N)$-random walks as $N \to \infty$, for sequences $(c_j)$ not depending on $N$. In this so-called {\em hierarchical mean field limit} (see \cite{[10]} and references therein), the reciprocal of the constant $\mu$ plays an important role as a scaling parameter concerning separation of time scales (Remark \ref{rem3.5.11}).
\end{remark}
\begin{example}\label{ex3.1.6}\rm
 The $(1, (c^j), N)$-random walk (called $c$-random walk in \cite{[9]}) has  jump probabilities
\begin{equation}\label{3.1.9}
r_j=\left(1-{c\over N}\right) \left({c\over N}\right)^{j-1},
\quad j=1,2,
\ldots,\quad{\rm where}\quad 0<c<N.
\end{equation}
In this case $h_j$ defined by (3.1.6) is given by
\begin{equation}\label{3.1.10}
h_j={N^2-c\over N(N-1)}\left({c\over N}\right)^{j-1},\quad j=1,2,\ldots .
\end{equation}
This random walk will sometimes be called {\it $c^j$-random} walk for brevity. Note also that by Remark \ref{rem3.1.5} a $(1, (c^j), N)$ random walk is the same as a $(\mu, (1), N)$ random walk with 
$$\mu = \frac{\log N}{\log(N/c)}.$$
\end{example}
\subsection {Degrees of hierarchical random walks}\label{sec3.2}
\setcounter{equation}{0}
Since $\Omega_N$ is countable and the random walks are
irreducible, it suffices to  consider,  instead of the operator $G^\zeta$ defined in (\ref{eq:2.1}), the number
$g^{(\zeta)}$ defined by (\ref{gzeta}).
 The following formula with discrete-time transition probabilities can also be used:
\begin{equation}\label{gzetaform}
g^{(\zeta)}=\frac 1{\Gamma(\zeta)}
\sum^\infty_{n=1}{\Gamma(\zeta+n)\over n!}p^{(n)}
(0,0).
\end{equation}
\begin{remark}\label{rem3.2.1}\rm
(a) We have from \cite{[9]} that
the $(1,(c^j),N)$-random walk has degree $\gamma^-$ with
\begin{equation}
\label{3.2.3}
\gamma ={\log c\over \log (N/c)},
\end{equation}
Equivalently (see Remark \ref{rem3.1.5}) the $(\mu, (1), N)$-random walk, $\mu > -1$, has degree $\gamma^-$ with $\gamma = \mu-1$. Note that the range of degrees of the $(1,(c^j),N)$-random walks is $(-1, \infty)$. In this sense this class is richer than the class of $\alpha$-stable processes on $\mathbb R^d$ (and $(\alpha,d)$-random walks).  

(b) Another consequence from \cite{[9]} is that for a $(1,(c^j),N)$-random walk with degree $\gamma$ we have
\begin{eqnarray}\label{asp}
 p_t(0,0) &\sim&  \mbox{ const }t^{-(\gamma +1)}h_t,\\ \nonumber
g_t^{(\gamma+1)}&\sim& \mbox{ const } \log t, \nonumber\\ \label{growthg}
g_t^{(\zeta)}&\asymp&   t^{\zeta-(\gamma+1)} \quad \mbox{ for } \zeta > \gamma +1,
\end{eqnarray}
where $h_t= h_t^{(\gamma)}$ is a slowly oscillating function (recall that $g_t^{(\zeta)}$ is defined by (\ref{defgtzeta})).\end{remark}
\begin{remark}\label{rem3.2.2}\rm
(a) Comparing (\ref{gg})  with (\ref{growthg}) we see that $(1, (c^j), N)$-random walks and $( \alpha,d)$-random walks with degree $\gamma$ (i.e. $c = N^{1-\alpha/d})$ have the same order of growth of $g_t^{(\zeta)}$ for $\zeta \ge \gamma +1$.

(b) The $(1,(c^j), N)$-random walks can also be compared to $\alpha$-stable processes in terms of the decay of the potential operators. For positive integer $k<\gamma+1$, the $k$-th power
$G_{N,\gamma}^k$ of the potential operator $G_{N,\gamma}$ of this hierarchical random walk has a kernel of the form
(see \cite{[9]}, (4.2.2))
$$G_{N,\gamma}^k(0,x)={\rm const} N^{-|x|(1-k/(\gamma+1))},$$
where $\gamma$ is the degree (\ref{3.2.3}).
If $\gamma=\frac{d}{\alpha}-1$ (hence $d>\alpha k)$, this can be written as
$$G_{N,\gamma}^k(0,x)={\rm const} \rho(x)^{-(d-\alpha k)},$$
where
$$\rho(x)=N^{|x|/d}.$$%
$\rho(x)$ is the ``Euclidean radial distance'' of $x$ from $0$, so that
 the volume of a ball of radius $\rho$ grows like $\rho^d$. Therefore the powers of the potential operator of the $(1,(c^j), N)$-random walk and the respective ones for the $\alpha$-stable process  have the same spatial asymptotic decay.

(c) For the $(\mu,(c^j),N)$-random walk with   $0< c<N^{1/\mu}$
the degree is 
\begin{equation}\label{cdeg}
\gamma = \gamma_N =
\frac{\mu-1+\mu \log c/\log N}{1-\mu \log c/\log N}.
\end{equation}
Hence 
$\gamma_{N}\rightarrow \mu-1\quad{\rm as}\quad N\rightarrow\infty$, more precisely, $\gamma_{N}=\mu-1+O(1/\log N)$ as $N\rightarrow \infty$.
In the case $\mu = 2$,
since the degree equals 1 for Brownian motion $(\alpha=2)$ on $\erre^4$ or
simple symmetric random walk on $\zet^4$ (see (\ref{eq:2.19})),
the  $N\rightarrow\infty$ limit behaviour of this hierarchical random walk
can be viewed as
 corresponding to Euclidean dimension $d=4$.
  This case plays a role in  the behaviour of 
two-level branching systems discussed in  \cite{[10]}.
Hierarchical models ``around dimension 4'' also play a prominent role in statistical physics \cite{[11]}.
\end{remark}

We turn now  to transience properties of the  $(\mu,(c_j),N)$-random walk.
We will sometimes write
$p_t, G^{\zeta}, G^{k}_t, D$ with a subscript or
 superscript $(\mu)$ when we need to emphasize the dependence on $\mu$.

We have, from (\ref{3.1.6}) and  (\ref{3.1.8}),
\begin{eqnarray}\label{3.2.7}
h_j&=&r_js_j, \quad \quad {\rm where}\\ \nonumber \quad s_j&=&{N\over N-1}+{1\over r_j}
\sum^\infty_{i=j+1}
r_i={N\over N-1}+{1\over c_{j-1}}\sum^\infty_{i=j+1}
{c_{i-1}\over N^{(i-j)/\mu}},\,
j=1,2,\ldots , \nonumber
\end{eqnarray}
therefore
\begin{equation}\label{3.2.8}
h_j= D{d_{j-1}\over N^{(j-1)/\mu}}
\quad{\rm where}\quad d_{j-1}=c_{j-1}s_j,\quad j=1,2,\ldots .
\end{equation}

We need conditions for finiteness of the powers
$G^{\zeta}$ in terms of the $h_j$.
\medskip
\begin{proposition}\label{prop3.2.3.} For any $\zeta> 0$,
\begin{eqnarray}\label{3.2.9}\mbox {(1)} \quad \quad \quad \quad 
G^\zeta=G^\zeta_{(\mu)}<\infty\quad \mbox{iff} \quad
\sum^\infty_{j=1}{1\over N^jh^\zeta_j}
<\infty\quad \mbox{ iff }\,\quad
\sum^\infty_{j=0}
{1\over N^{j(\mu-\zeta)/\mu}d^\zeta_{j}}<\infty,
\end{eqnarray}
where $h_j$ and $d_j$ are given by (3.1.6), (3.2.7) and (3.2.8).

(2) In terms of the
 $c_j$ in  (3.1.8),
\begin{equation}\label{3.2.10}
G^\zeta=
G^\zeta_{(\mu)}<\infty\quad\,{\it iff}\quad\,\sum^\infty_{j=0}
{1\over N^{j(\mu-\zeta)/\mu}c^\zeta_{j}}
<\infty,
\end{equation}
 provided that
\begin{equation}\label{3.2.11}
\limsup_{j\rightarrow\infty}{1\over r_j}\sum^\infty_{i=j+1}r_i<\infty \quad\,
(or \, \,equivalently \quad\,
\limsup_{j\rightarrow\infty}{1\over c_{j}}
\sum^\infty_{i=j+1}{c_{i}\over N^{(i-j)/\mu}}<\infty).
\end{equation}
A sufficient condition for (\ref{3.2.11}) is
\begin{equation}\label{3.2.12}
\limsup_{j\rightarrow\infty}{c_{j+1}\over c_j}<N^{1/\mu},
\end{equation}
and hence for large $N$ it suffices that
$\limsup\limits_jc_{j+1}/ c_j<\infty$.
\end{proposition}
\begin{proof} (1)
We have,  from (\ref{3.1.5}) and (\ref{eq:2.1}),
\begin{equation}\label{3.2.13}
G^\zeta=G^\zeta_{(\mu)}=(N-1)\sum^\infty_{j=1}{1\over N^jh^\zeta_j},
\end{equation}
Then the first part of (\ref{3.2.9}) is obvious from (\ref{3.2.13}), and the second one follows from (\ref{3.2.8}).

\noindent
(2) (\ref{3.2.10}) and (\ref{3.2.11}) follow from part (1) and (\ref{3.2.7}). 
\end{proof}

The defining quantities of the $(\mu,(c_j),N)$-random walk  (\ref{3.1.8}) are $\mu,$ the sequence $(c_j)$ and $N$, but conditions for finiteness of the powers of $G$ are more conveniently established  by the sequence $(d_j)$ 
defined  by (\ref{3.2.7}) and (\ref{3.2.8}).
The $c_j$ can be obtained from the $d_j$ by (\ref{3.1.7}), (\ref{3.1.8}) and (\ref{3.2.8}); e.g., for $\mu=1$,
\begin{equation}\label{3.2.14}
c_j=\frac{N-1}{N}d_j-\frac{(N-1)^2}{N^3}N^{2j}\sum^\infty_{i=j}\frac{d_i}{N^{2i}}.
\end{equation}

An obvious consequence of the previous proposition for the
$(\mu,(c_j),N)$-random walk is
\begin{equation}\label{3.2.15'}
G^\mu < \infty \mbox { iff}\quad\sum^\infty_{j=1}{1\over d^{\mu}_j}<\infty ,
\end{equation}
or in terms of the $c_j$,
\begin{equation}\label{3.2.16'}
G^\mu < \infty \mbox { iff}\quad\sum^\infty_{j=1}{1\over
c^{\mu}_j}<\infty,
\end{equation}
{\rm provided\,\,that}
\begin{equation}\label{3.2.17}
\limsup_{j\rightarrow\infty}{1\over c_j}\sum^\infty_{i=j+1}{c_i\over N^{(i-j)/\mu}}<\infty.
\end{equation}
Note that (\ref{3.2.7}) and (\ref{3.2.8}) imply $d_j\geq c_j$, therefore $\sum_j 1/c^\mu_j<\infty$ implies $G^\mu < \infty$.

%The $(1,(c^j),N)$-random walk trivially satisfies condition (3.2.12), and represented as a \linebreak $(k+1,(c_j),N)$-random walk, from (3.1.9), (3.1.10) and (3.2.8) we have
%\begin{equation}
%d_j=\frac{N^2-c}{(N-1)(N-c)}\left(\frac{c}{N^{k/(k+1)}}\right)^j.
%\end{equation}

We have seen that 
 the  $c^j$-random walk with degree $\gamma$
actually has degree $\gamma^-$. Now we ask for existence of
$(1,(c_j),N)$-random walks of degree $\gamma^+$.
The next proposition and its corollary give an answer.
\vglue.45cm
\noindent
\begin{proposition}\label{prop3.2.4}
Consider a $(\mu,(c_j),N)$-random walk  such that 
\begin{equation}\label{3.2.19}
\inf_jc_j>0\quad{\it and}\quad
\lim_{j\rightarrow\infty}\frac{c_{j+1}}{c_j}=1.
\end{equation}
Then for each $\varepsilon > 0$
\begin{equation}\label{nearmuupper}
p_t(0,0) = o(t^{-\mu+\varepsilon}) \mbox{ as } t \to \infty
\end{equation}
and 
\begin{equation}\label{nearmulower}
p_t(0,0)^{-1} = o(t^{\mu+\varepsilon}) \mbox{ as } t \to \infty.
\end{equation}
\end{proposition}
\begin{proof}
It is not difficult to show from (\ref{3.2.7}) and (\ref{3.2.8}) that (\ref{3.2.19}) implies
\begin{equation}
\label{Eq:3.2.20}
\lim_{j\rightarrow\infty}\frac{d_{j+1}}{d_j}=1.
\end{equation}
We have from (\ref{3.1.5}) and (\ref{3.2.8}) that 
\begin{eqnarray*}
p_t(0,0) &=& \mathrm {const} \sum_{j=1}^\infty N^{-j} \exp\left(-D\frac{d_{j-1}}{N^{(j-1)/\mu}}t\right)\\
&\le & \mathrm {const} \sum_{j=1}^\infty N^{-j} \exp\left(-K\frac{c^{j-1}}{N^{(j-1)/\mu}}t\right),
\end{eqnarray*}
where $c\, (< 1)$ can be chosen arbitrarily close to 1 and $K$ is another constant. The latter expression can be rewritten as 
$$\mathrm {const} \sum_{j=1}^\infty N^{-j} \exp\left(-K\frac{1}{N^{(j-1)/\mu'}}t\right),$$
where $\mu' \,(< \mu)$ can be chosen arbitrarily close to $\mu$. The estimate (\ref{nearmuupper}) is now immediate from (\ref{asp}), where the constant $\gamma$ appearing there is chosen as $\mu'-1$. The estimate (\ref {nearmulower}) is proved in an analogous way.
\end{proof}

\begin{corollary}\label{cor3.2.5}
Under the assumptions of Proposition \ref{prop3.2.4},

(a) the degree of the random walk  is $\gamma=\mu-1$, and it is 
$\gamma^+$  iff $\sum\limits_j1/c^\mu_j<\infty$, 

(b) if $\mu < 1$ then the return time to $0$ has all moments of order less than $1-\mu$,

(c) if $\mu = 1$ then the return time to $0$ has no moments of positive order.
\end{corollary}
\medskip
\begin{proof} (a) is immediate from (\ref{nearmuupper}), (\ref {nearmulower}), (\ref{defdegree}) and by noting that condition (\ref{3.2.12}) holds.

\noindent
(b) follows from (\ref {nearmulower}) and Proposition \ref{prop2.2.5}.

\noindent
(c) follows from (\ref{nearmuupper}) and Proposition \ref{nomoments}.
\end{proof}
\medskip

Next we give an example of a walk satisfying the assumptions of Proposition \ref{prop3.2.4}. 
\begin{example}\label{ex3.2.6}\rm
($j^\beta$-{\it random walk}). 
Consider a $(\mu ,(c_j),N)$-random walk with $\mu > 0$ such that $d_j$ defined by (\ref{3.2.7}) and (\ref{3.2.8}) is given by
\begin{equation}\label{3.2.21}
d_j=(j+1)^\beta,\quad j=0,1,\ldots, \quad{\rm where}\quad \beta\geq 0.
\end{equation}
We call this a $j^\beta$-{\it random walk}, referring to $d_j$ rather than to 
$c_j$. The degree of this random walk is $\gamma = \mu-1$, and it is $\gamma^+$ if $\beta > 1/\mu$, and $\gamma^-$ if $\beta \le 1/\mu$.
 Since $c_j$ also behaves like $(j+1)^\beta$ (see (\ref{djcj}) below), we could consider the random walk with $c_j=(j+1)^\beta$ instead of (\ref{3.2.21}), but this would complicate the exact asymptotics  derived in subsection \ref{sec3.3} because they are more directly related to $d_j$ than to $c_j$.
\end{example}
 
 The next result shows in particular that for a large class of sequences $(c_j)$ the degree of the $(\mu,(c_j),N)$-random walk approaches $\mu$ as $N\to \infty$.
 \medskip 
  \begin{proposition}\label{prop3.2.8} Consider a $(\mu,(c_j),N)$-random walk, denote its degree by $\gamma$ and put
$$\overline c=\limsup_{j\rightarrow\infty}\frac{c_{j+1}}{c_j}, \quad \underline c = \liminf_{j\rightarrow\infty}\frac{c_{j+1}}{c_j}. $$ 

(1) If $\overline c < N^{1/\mu}$, then
\begin{equation}\label{3.2.22}\gamma \le \frac{\mu-1+\mu\log \overline c/\log N}{1-\mu\log \overline c/\log N}.
\end{equation} 

(2) If $\underline c < N^{1/\mu}$, then
\begin{equation}\label{lowerc}\gamma \ge \frac{\mu-1+\mu\log \underline c/\log N}{1-\mu\log \underline c/\log N}.
\end{equation}

(3) If $0 < \underline c \le \overline c < \infty$, then
$\gamma =\mu-1+O(1/\log N)$ as $N\rightarrow\infty$.
\end{proposition}
\begin{proof} (1) For each $a \in (\overline c, N^{1/\mu})$, the $(\mu,(a^j),N)$-random walk has degree (see (\ref{cdeg}))
$$\gamma^{(a)}=\frac{\mu-1+\mu\log a/\log N}{1-\mu\log a/\log N}.$$
Let $0<\zeta<\gamma+1$. Then, by (\ref{eq:2.1}), (\ref{3.2.10}) and by the definition of 
the degree $\gamma$,
\begin{equation} \label{finitesum}
\sum_j\frac{1}{N^{j(\mu-\zeta)/\mu}c^\zeta_j}<\infty.
\end{equation}
Since  $c_j \le Ka^j$ for all $j=0,1,...$  and a suitable constant $K> 0$, (\ref{finitesum}) implies
$$\sum_j\frac{1}{N^{j(\mu-\zeta)/\mu}a^{j\zeta}}<\infty.$$
Consequently, $\zeta\leq \gamma^{(a)}+1$. It follows that $\gamma\leq \gamma^{(a)}$, and since $a$ is arbitrary, the assertion (\ref{3.2.22}) follows.

(2) is proved in an analogous way, and (3) is immediate from (1) and (2).
\end{proof}

We now pass to last exit times. The following results describe
the behaviour of  moments of the last exit time $L_{B_R}$ from a closed  ball $B_R$  of radius $R$ for transient $c^j$-random walks.
\medskip
\begin{proposition}\label{prop3.2.10} For a $(\mu,(\eta^j),N)$-random walk 
with $\mu \ge 1$, $\eta > 1$ and  $B_R$ a closed ball of  radius $R$ centered at $0$,
\begin{equation}\label{3.2.25}
\int^\infty_0t^{\mu-1}P_t(0,B_R)dt
=\Gamma(\mu){(N-1)^{\mu+1}\over (N^{(\mu+1)/\mu}\eta^{-1}-1)^{\mu}}{1\over \eta^\mu-1}
\left({N\over \eta^\mu}\right)^{R}.
\end{equation}
where $P_t(0,B_R):=\sum_{x\in B_R}p_t(0,x)$.
\end{proposition}
\begin{corollary} \label{cor3.2.11}  Under the conditions of Proposition \ref{prop3.2.10},\\
(1)  for fixed $N$,
\begin{equation}\label{3.2.26}
\ee_0L^{\mu-1}_{B_R}\asymp
\Gamma(\mu){(N-1)^{\mu+1}\over (N^{(\mu+1)/\mu}\eta^{-1}-1)^{\mu}}{1\over \eta^\mu-1}\left({N\over \eta^\mu}\right)^{R}
\quad\mbox{\it as}\quad R\rightarrow\infty,
\end{equation}
 (2)  
For fixed $R$,
\begin{equation}\label{3.2.27}
\ee_0L^{\mu-1}_{B_R}\asymp\Gamma(\mu){\eta^{\mu}\over (\eta^{\mu}-1)}\left(\frac{N}{\eta^{\mu}}\right)^R
\quad\mbox{\it as}\quad N\rightarrow\infty,
\end{equation}
and in particular \\
(3)  
\begin{equation}\label{3.2.28}
\frac{\ee_0L^{\mu-1}_{B_{R+1}}}{\ee_0L^{\mu-1}_{B_R}}\asymp\frac{N}{\eta^{\mu}}\quad
{ as}\quad N\rightarrow\infty.
\end{equation}
\end{corollary}
{\it Proof of Proposition \ref{prop3.2.10} and Corollary \ref{cor3.2.11}:}
We sketch the proof of (\ref{3.2.25}).

Writing $c = \eta N^{(\mu-1)/\mu}$,  $h_j=ba^{j-1}$ with $b=(N^2-c)/N(N-1)$ and $a=c/N$ (see (\ref{3.1.10})) we obtain from (\ref{3.1.5})
\begin{eqnarray*}
\lefteqn{
\int^\infty_0t^{\mu-1}P_t(0,B_R)dt=\sum_{x\in B_R}\int^\infty_0t^{\mu-1}p_t(0,x)dt}\\
&=&\frac{\Gamma(\mu)}{b^{\mu}}
\biggl[
(N-1)\sum^\infty_{j=1}\frac{1}{N^ja^{(j-1)\mu}}-\sum^R_{m=1}\frac{N^{m-1}(N-1)}{N^ma^{(m-1)\mu}}\biggr.\\
&&+(N-1)\sum^R_{m=1}N^{m-1}(N-1)\sum^\infty_{j=m+1}\biggl.\frac{1}{N^ja^{(j-1)\mu}}
\biggr].
\end{eqnarray*}
Computing the summations and substituting the expressions for $a$ and $b$ leads to (\ref{3.2.25}).

The results (\ref{3.2.26}) and (\ref{3.2.27}) are obtained from (\ref{3.2.25}) and
(\ref{eq:2.5}), and
(\ref{3.2.28}) also follows  from the previous results.
\hfill$\Box$
\begin{remark}\label{rem3.2.12} \rm
(a) Consider a $c^j$-random walk on $\Omega_N$ and an $\alpha$-stable process on $\mathbb R^d$ having the same degree $\gamma$. 
We see from the previous corollary that, for $\zeta < \gamma$,
$\ee_0L^\zeta_{B_R}$ grows like $R^{\alpha \zeta}$
for the $\alpha$-stable process on $\erre^d$,
 and it  grows like $(N/c)^{(\zeta+1)R}$ for the $c^j$-random walk on $\Omega_N$.
If the degrees of the two processes  coincide (i.e., $c=N^{1-\alpha / d}$) then $(N/c)^{\zeta R}=\rho^{\alpha \zeta}$ where $\rho=N^{R/d}$ is the Euclidean radial distance from 0 (Remark \ref{rem3.2.2}(b)). 
This shows that a $c^j$-random walk takes on the average a longer time to leave an ``Euclidean'' ball than an $\alpha$-stable process with the same degree.

\noindent
(b) Part (3) of Corollary \ref{cor3.2.11}  shows that there is separation of time scales on balls of hierarchical radius $R$ and $R+1$ as $N\rightarrow\infty$. The analogue to (\ref{3.2.28})  is true for the $\alpha$-stable process on $\erre^d$, on balls of ``Euclidean'' radius $N^{R/d}$ and $N^{(R+1)/d}$ (corresponding to hierarchical distance $R$ and $R+1$) as $N\rightarrow\infty$. Indeed, on appropriate time scales, one sees certain features of Euclidean random walks
related to separation of time scales. See e.g. Cox and Griffeath \cite{CG86} where diffusive clustering in the two-dimensional voter model is shown based on such features of two-dimensional simple random walk.

%

%On the other hand, for the $\alpha$-stable process on $\erre^d$ the left-hand side of (\ref{3.2.28}) tends to 1 as $R\rightarrow\infty$, and there is no separation of time scales. 

\noindent
(c) For the $c^j$-random walk, from (\ref{3.1.9}) we have  $r_j=(N/c)r_{j+1}$. Hence a jump of size $j$  is $N/c$ times more likely than a jump of size $j+1$. Therefore, as time flows the points visited by the random walk form a clustered pattern:  the  walk spends some (long) time jumping within a closed ball of radius
$j$, and forming a cluster there, before jumping to a point outside the ball, and beginning a new cluster within another ball of radius $j$, which by the ultrametric structure of $\Omega_N$
 is necessarily disjoint from the  previous ball, and so on. This behaviour is analogous to that of the Weierstrass random walk on the lattice studied in \cite{[30], [31], [32]}. The one-dimensional Weierstrass random walk has step distribution with density function
 $$\frac{a-1}{2a} \sum_{n=0}^\infty a^{-n}[\delta(x-\Delta b^n) + \delta(x+\Delta b^n)], \quad x \in \mathbb R,$$
 where $a$, $b$ and $\Delta$ are constants, $a > 1$, $b > 0$, $\Delta > 0$. When $b$ is an integer the walk stays on a lattice. (The characteristic function of the step distribution is Weierstrass' example of a function which is everywhere continuous and nowhere differentiable.) The $d$-dimensional Weierstrass random walk is an obvious extension.
\end{remark}
\subsection {A special   class of hierarchical random walks} \label{sec3.3}
\setcounter{equation}{0}
We know that $(\mu,(1), N)$-random walks have degree $(\mu-1)^-$
and  $g^{(\mu)}_t$ defined by (\ref{defgtzeta}) grows logarithmically.
 In this subsection we will construct a class of hierarchical walks with degree  $(\mu-1)^-$ for which $g^{(\mu)}_t$ grows only {\em sublogarithmically}.
To this end we consider $(\mu,(c_j), N)$-random walks defined by (\ref{3.1.8}) 
such that $c_j\leq c_{j+1}$ for all $j$. It can be shown
easily from (\ref{3.2.7}) and (\ref{3.2.8}) that this assumption implies $d_j \le  d_{j+1}$ 
for all $j$ and
\begin{equation}\label{djcj}
     \frac{N}{N-1} < \frac{d_j}{c_j} < \frac{N}{N-1} + \frac1{N^{1/\mu-1}} \quad \mbox { for all } j.
     \end{equation}

If $c_j$ is non-decreasing, then $\liminf c_{j+1}/c_j \ge 1$. Hence Proposition \ref {prop3.2.8} implies that the degree is greater or equal to $\mu-1$. If we assume in addition that $\sum_j  1/d_j^{\mu} = \infty$, then  (\ref{3.2.15'}) implies that $G^\mu$ is infinite, hence the degree is $(\mu-1)^-$.

To state the next proposition in a compact way, we put
\begin{equation}\label{defH}
f^{(1)}_t = G_t(0,0), \quad f^{(2)}_t = G_t^2(0,0), \quad
f^{(3)}_t = (G_t^2G)(0,0).
\end{equation}
\begin{proposition}\label{prop3.3.1} Assume  $d_j \leq d_{j+1}$ for all $j$, and
\begin{equation}\label{recmu}
\sum_j \frac 1{d_j^{\mu}} = \infty.
\end{equation}
(a) In case $\mu =1, 2$ or $3$,
$$f^{(\mu)}_t\sim\,\displaystyle{N-1\over ND_{(\mu)}^\mu}\sum\limits^{\mu \log t/\log N}_
{j=0}\displaystyle{1\over {d_j^\mu}}\quad as \quad t\rightarrow\infty$$
where $D_{(\mu)}$ is the normalizing constant  in (\ref{3.1.8}).
\\(b) For general $\mu$ and $g_t^{(\mu)}$ defined by (\ref{defgtzeta}),
$$g_t^{(\mu)} \sim \displaystyle{N-1\over ND_{(\mu)}^\mu} \sum\limits^{\mu \log t/\log N}_
{j=0}\displaystyle{1\over {d_j^\mu}}\quad as \quad t\rightarrow\infty.$$
 (The upper limits in the sums are understood as integer part.)
 \end{proposition}
\medskip
\begin{proof} 
Denote (see (\ref{3.1.5}) and (\ref{3.2.8}))
\begin{equation}\label{3.3.1}
p^{(\mu)}_t=p^{(\mu)}_t(0,0)={N-1\over N}q^{(\mu)}_t,
\end{equation}
\begin{equation}\label{3.3.2}
q^{(\mu)}_t=\sum\limits^\infty_{j=0}
{{\rm exp}
 \{-
  {d_j\over N^
   {j/\mu}}D_{(\mu)} t\}
  \over N^j
 },
\end{equation}
We will omit the superscript and subscript $(\mu)$  but the value of 
$\mu$ will be clear in each case.

(a) {\em Case} $\mu=1$.  By (\ref{3.3.1}), (\ref{3.3.2}) (\ref{defH}) and (\ref{Gt}),
$$
G_t(0,0)={N-1\over ND}\int^{Dt}_0q_sds,\quad \hbox{\rm where} \quad q_t=\sum_j
{{\rm exp}\{-{d_j\over N^j}t\}\over N^j}.$$

The Laplace transform of $q_t$ is
$$
\widetilde{q}(\lambda)=\int^\infty_0e^{-\lambda t}q_tdt=
\sum_j{1\over N^j}{1\over (\lambda+{d_j\over N^j})}=\sum_j{1\over \lambda N^j+d_j},
$$
and $\widetilde{q}(\lambda)\rightarrow\infty$ as $\lambda\rightarrow 0$ by
(3.3.3).

Write
$$
\widetilde{q}(\lambda)=F_1(\lambda)+F_2(\lambda),
$$
where
$$
F_1(\lambda)=\sum_{j\leq Q(\lambda)}{1\over \lambda N^j+d_j},\quad F_2(\lambda)=\sum_{j>Q(\lambda)}{1\over \lambda N^j+d_j},$$
with
$$Q(\lambda)=-{\log \lambda\over \log N},\quad 0<\lambda < 1.$$
Since
$$F_2(\lambda)\leq {1\over \lambda}\sum_{j>Q(\lambda)}{1\over N^j}\leq
L {1\over \lambda N^{Q(\lambda)}}= L,$$
where $L$ is a constant, then
$$\widetilde{q}(\lambda)\sim F_1(\lambda)\quad{\rm as}\quad \lambda\rightarrow 0.$$%

Write
$$F_1(\lambda)=J_1(\lambda)+J_2(\lambda),$$
where
\begin{eqnarray*}
J_1(\lambda)&=&\sum_{j\leq Q(\lambda)}{1\over d_j},\\
J_2(\lambda)&=&\sum_{j\leq Q(\lambda)}\left({1\over \lambda N^j+d_j}-{1\over d_j}\right)=-\sum_{j\leq Q(\lambda)}{\lambda N^j\over (\lambda N^j+d_j)d_j}.
\end{eqnarray*}
Since $\inf_j d_j>0$,
$$|J_2(\lambda)|\leq L \lambda\sum_{j\leq Q(\lambda)}N^j\leq
L_1 \lambda N^{Q(\lambda)}=L_1,$$
where $L$ and $L_1$ are constants, then
$F_1(\lambda)\sim J_1(\lambda)$, and therefore
$$\widetilde{q}(\lambda)\sim J_1(\lambda)\quad{\rm as}\quad \lambda\rightarrow 0.$$

Let
$$H(t)=\sum_{j\leq Q(t^{-1})}{1\over d_j},\quad t>0,$$
so $J_1(\lambda)=H(1/\lambda)$. $H(t)$ is slowly varying at $\infty$. Indeed, let $x>1$,
 then
$${H(tx)\over H(t)}=1+\sum_{j}R_{t,x}(j),$$
where
$$R_{t,x}(j)=
{{d^{-1}_j}\UNO{[Q(t^{-1})<j\leq Q((tx)^{-1})]}\over
\sum_{k\leq Q(t^{-1})} d_k^{-1}}.$$%
Since the sequence $d_j$ is non-decreasing,
$$\sum_j R_{t,x}(j)\leq
{
d^{-1}_{Q(t^{-1})+1}(Q((tx)^{-1})-Q(t^{-1}))\over
 d^{-1}_{Q(t^{-1})}Q({t^{-1})}
}
\leq{\log(tx)-\log t\over \log t}={\log x\over \log t}
\rightarrow 0\quad{\rm as}\quad t\rightarrow\infty,$$
hence
$${H(tx)\over H(t)}\rightarrow 1\quad{\rm as}\quad t\rightarrow\infty.$$
A similar argument works for $0<x<1$.

By a Tauberian theorem (\cite{[0]}, Theorem 1.7.1)
$$\int^t_0 q_s ds \sim
\sum^{Q(t^{-1})}_{j=1}{1\over d_j}\quad{\rm as}\quad t\rightarrow\infty,$$%
and the conclusion follows.

{\em Case} $\mu = 2$.  Using  the formula
$$G^{2}_t(0,0) =\int^t_0
\int^t_0p_{s+r}dsdr=2
\int^t_0
\int^r_0p_{s+r}dsdr
$$
we have
$$
G^2_t(0,0) =2{N-1\over ND^2}
\int^{Dt}_0\int^r_0q_{s+r}dsdr,\quad {\rm where} \quad  q_t=\sum_{j}
{{\rm exp}\{-{d_j\over N^{j/2}}t\}\over N^j}.
$$

Let
\begin{eqnarray*}
M_t=\int^t_0q_{s+t}ds=\int^{2t}_tq_sds&=&\sum_j{1\over N^j}
\int^{2t}_t
{\rm exp}\left\{-
   {d_j\over N^{j/2}
   }s\right\}
  ds\\
&=&\sum_j
{{\rm exp}
 \{-
  {d_j\over N^{j/2}
  }t\}
-{\rm exp}
 \{-
  {d_j\over N^
   {j/2}
  }2t\}
 \over
  N^{j/2}d_j}.
\end{eqnarray*}

The Laplace transform of $M_t$ is
\begin{eqnarray*}
\widetilde{M}(\lambda)&=&\sum_j
{1\over N^{j/2}d_j}
\left({1\over \lambda+{d_j\over N^{j/2}}}-{1\over \lambda +
2{d_j\over N^{j/2}}}\right)
\\
&=&\sum_j{1\over (\lambda N^{j/2}+d_j)(\lambda N^{j/2}+2d_j)} 
=F_1(\lambda)+F_2(\lambda),
\end{eqnarray*}
where
\begin{eqnarray*}
F_1(\lambda)&=&\sum_{j\leq Q(\lambda)}{1\over (\lambda N^{j/2}+d_j)(\lambda
N^{j/2}+2d_j)},\\
F_2(\lambda)&=&\sum_{j>Q(\lambda)}{1\over (\lambda N^{j/2}+d_j)(\lambda N^{j/2}+2d_j)},
\end{eqnarray*}
with
$$Q(\lambda)=-2{\log\lambda\over \log N},\quad 0<\lambda<1.$$%
Since
$$F_2(\lambda)\leq {1\over\lambda^2}\sum_{j>Q(\lambda)}{1\over N^j}\leq
L {1\over \lambda^2N^{Q(\lambda)}}=L,$$
then
$$\widetilde{M}(\lambda)\sim F_1(\lambda)\quad{\rm as}\quad \lambda\rightarrow 0.$$%

Write
$$
F_1(\lambda)=J_1(\lambda)+J_2(\lambda),$$
where
$$J_1(\lambda)={1\over 2}\sum_{j\leq Q(\lambda)}{1\over d^2_j},$$
\begin{eqnarray*}
J_2(\lambda) &=&\sum_{j\leq Q(\lambda)}\biggl(
{
1\over (\lambda N^{j/2}+d_j)(\lambda N^{j/2}+2d_j)
} -{1\over 2d^2_j}\biggr)\\
&=&- \sum_{j\leq Q(\lambda)}
{\lambda^2N^{j}+3d_j \lambda N^{j/2}\over 2(\lambda N^{j/2}+d_j)(\lambda N^{j/2}+2d_j)d^2_j}.
\end{eqnarray*}
Since
$$
|J_2(\lambda)|\leq L\!
\biggl(\lambda^2\sum_{j\leq Q(\lambda)}N^{j}+\lambda
\sum_{j\leq Q(\lambda)}N^{j/2}\biggr)
\leq
L_1 (\lambda^2N^{Q(\lambda)}+\lambda
N^{Q(\lambda)/2}
)=L_2,
$$
then
$
F_1(\lambda)\sim J_1(\lambda)$  and therefore
$$\widetilde{M}(\lambda)\sim J_1(\lambda)\quad{\rm as}\quad \lambda\rightarrow 0.
$$

Let
$$
H(t)={1\over 2}\sum_{j\leq Q(\lambda)}{1\over d^2_j},
$$
so $J_1(\lambda)=H(1/\lambda)$. $H(t)$ is slowly varying at $\infty$ (as above), and the conclusion follows by follows by the
Tauberian theorem.

The case $\mu =3$ is proved similarly, using the formula
$$G^2_tG=2\int^t_0\int^r_0\int^\infty_0p_{s+r+u}dudsdr.$$

(b) The proof is analogous to that of part (a) for $\mu=1$, hence we will only give a sketch showing a step which is different.

Proceeding as above we obtain from (\ref{defgtzeta}), (\ref{3.3.1}) and (\ref{3.3.2}),
$$g^{(\mu)}_t=\frac{1}{\Gamma(\mu)}\frac{N-1}{N}\frac{1}{D^\mu}
\int^{Dt}_0q_sds ,$$
where
$$q_t=t^{\mu-1}\sum_j\frac{{\rm exp}\{-\frac{d_j}{N^{j/\mu}}t\}}{N^j}.$$

The Laplace transform of $q_t$ is given by
$$\widetilde{q}(\lambda)=\Gamma(\mu)\sum_j\frac{1}{(\lambda N^{j/\mu}+d_j)^\mu}=\Gamma(\mu)(F_1(\lambda)+F_2(\lambda)),$$
where
$$F_1(\lambda)=\sum_{j\leq Q(\lambda)}\frac{1}{(\lambda N^{j/\mu}+d_j)^\mu},\quad F_2(\lambda)=\sum_{j>Q(\lambda)}\frac{1}{(\lambda N^{j/\mu}+d_j)^\mu},$$
with
$$Q(\lambda)=-\mu\frac{\log \lambda}{\log N},\quad 0<\lambda <1,$$
and $F_2(\lambda)$ is bounded, so $q(\lambda)\sim \Gamma(\mu)F_1(\lambda)$ as $\lambda\rightarrow 0$.

Write
$$F_1(\lambda)=J_1(\lambda)+J_2(\lambda),$$
where
$$J_1(\lambda)=\sum_{j\leq Q(\lambda)}\frac{1}{d^\mu_j},\quad
J_2(\lambda)=\sum_{j\leq Q(\lambda)}
\frac{d^\mu_j-(\lambda N^{j/\mu}+d_j)^\mu}{
(\lambda N^{j/\mu}+d_j)^\mu d^\mu_j}.$$

We show that $J_2(\lambda)$ is bounded.

\noindent
{\it Case $\mu>1$}: By convexity,
$$(a+b)^\mu-b^\mu\leq 2^{\mu-1}a^\mu+(2^{\mu-1}-1)b^\mu,\quad a,b\geq 0.$$
Using this inequality with $a=\lambda N^{j/\mu}$ and $b=d_j$ we obtain
\begin{equation}\label{J1}
|J_2(\lambda)|\leq\sum_{j\leq Q(\lambda)}\left[\frac{2^{\mu-1}\lambda^\mu N^j}{(\lambda N^{j/\mu}+d_j)^\mu d^\mu_j}+
\frac{(2^{\mu-1}-1)}{(\lambda N^{j/\mu}+d_j)^\mu}\right].
\end{equation}
{\it Case $0<\mu<1$:} Using the obvious inequality
$$(a+b)^\mu-b^\mu\leq a^\mu,\quad a,b\geq 0,$$
with $a=\lambda N^{j/\mu}$ and $b=d_j$ we obtain
\begin{equation}\label{J2}
|J_2(\lambda)|\leq \sum_{j\leq Q(\lambda)}\frac{\lambda^\mu N^j}
{(\lambda N^{j/\mu}+d_j)^\mu d^\mu_j}.
\end{equation}
Inequalities (\ref{J1}) and (\ref{J2}) imply that $J_2(\lambda)$ is bounded in both cases.

Therefore $\tilde{q}(\lambda)\sim \Gamma(\mu)J_1(\lambda)$ as $\lambda\rightarrow 0$, and the rest of the proof is like that of part (a) for $\mu=1$.
\end{proof}
\begin{remark}\label{rem3.3.3}\rm
(a) In \cite{[9]} we derived exact asymptotics for the growth of the incomplete potential operators $G_t$ for recurrent $c^j$-random walks. Unless $c=1$, these walks  have degree $< 0$ (see (\ref{3.2.3})) and hence  behave differently from the critically recurrent
walks inverstigated in Proposition \ref{prop3.3.1} (case $\mu=1$).\\
\noindent (b)
The proof of Proposition \ref{prop3.3.1} for $\mu=1$ provides a form of approximation for a class of divergent series,  including the series $\sum n^{-s}$, $0<s\leq 1$, related to the Riemann Zeta function \cite{[25]}.
\end{remark}

Using the well known formulas
$$\sum^n_{j=1}{1\over j}\sim \log n
\quad{\rm and}\quad
\sum^n_{j=1}{1\over j^\beta}\sim {n^{1-\beta}\over 1-\beta}\quad{\rm for}\,\,
\beta\in(0,1)\quad{\rm as}\quad n\rightarrow\infty,$$
\noindent
we obtain the following results from Proposition \ref{prop3.3.1}:

\noindent
\begin{corollary}\label{cor3.3.4} The
$(\mu, ((j+1)^\beta), N)$-random walk (with $0<\beta$) has degree
$\gamma = \mu-1$, and it has degree $\gamma^-$ iff $\beta \le
\mu^{-1}$. In this case, $g_t^{(\mu)}$ grows like const $\log \log t$
for $\beta = \mu^{-1}$, and like const $(\log t)^{1-\beta \mu}$
for $0 < \beta < \mu^{-1}$. Note that these growths have a similar
pattern as (\ref{eq:2.22}) and (\ref{gg}) for the 
 $\alpha$-stable process and the
$(\alpha,d)$-random walk, and
(\ref{growthg}) for the $c^j$-random walk, except that $t$ is now replaced by $\log t$.
\end{corollary}

The $j^\beta$-random walk defined in Example \ref{ex3.2.6} is a special case for Proposition \ref{prop3.3.1}, and we obtain from it as an ingredient for our discussion of occupation time fluctuations 
 of $j^\beta$-branching random walks (subsection \ref{sec3.6}) the following exact asymptotics: 
\begin{equation}\label{3.3.12}
\begin{array}{rcllll}
\mu =1,2:
&G^{\mu}
&<\infty
&\mbox{\rm for }\beta>\frac{1}{\mu},\\
&G^{\mu}_t
&\sim\frac{N-1}{ND^{\mu}}\log\log t
&\mbox{\rm for }
\beta=\frac{1}{\mu},\\
&G^{\mu}_t
&\sim
\frac{(N-1)\mu^{1-\mu\beta}}{ND^{\mu}(1-\mu\beta)(\log N)^{(1-\mu\beta)}}(\log t)^{1-\mu\beta}
&\mbox{\rm for }
0<\beta<\frac{1}{\mu},\\
\mu = 3:&G^3
&<\quad\infty\quad
&\mbox{\rm for } \beta > {1\over 3},\nonumber\\
&
G^{2}_tG
&\sim\quad{N-1\over ND^{3}}\log\log t\quad
&\mbox{\rm for } \beta={1\over 3},\nonumber\\
&
G^{2}_tG
&\sim\quad{(N-1)3^{1-3\beta}\over ND^{3}(1-3\beta)(\log
N)^{1-3\beta}}(\log t)^{1-3\beta}
&\mbox{\rm for } 0<\beta<{1\over
3}.
\end{array}
\end{equation}
(Recall that $D^\mu=(D_{(\mu)}^\mu,\mu=1,2,3$).
   
\subsection {An occupation time limit}\label{sec3.4}
\setcounter{equation}{0}
The incomplete potential operator $G_t$ defined by (\ref{Gt}) is also the norming for occupation time limits of Darling-Kac type \cite{[8]} for recurrent random walks. For the 
critically recurrent random walks of  subsection 3.3 we have the following result:

\begin{proposition}\label{prop3.4.1} Let $X=\{X(t), t\geq 0\}$ be the continuous time version of the $(1, (c_j), N)$-random walk with $c_j \leq c_{j+1}$ for all $j$ in the recurrent case ($\sum_j d^{-1}_j=\infty$ where $d_j$ is given by (\ref{3.2.7}), (\ref{3.2.8})). Then for any function $F:\Omega _N\rightarrow \erre^+$ with bounded support,
\begin{equation}
\label{eq:3.4.1}
\pee \left[ \frac{ND}{(N-1)\sum_{y\in \Omega _N}F(y) \sum_{j\leq \log t/\log N}d^{-1}_j}\int^t_0 F(Y(s))ds<x\right]\rightarrow 1-e^{-x},\;\; x\geq 0,
\end{equation}
as $t\rightarrow \infty$, where $D$ is the normalizing constant in (\ref{3.1.8}). 
\end{proposition}
\begin{proof} Using (\ref{3.1.5}) we have for $\lambda >0$,
\begin{eqnarray*}
\pi_{\lambda}(x,y)&:=& \int^\infty_0 e^{-\lambda t}p_t (x,y)dt\\
&=& (\delta_{0, |x-y|}-1 )\frac{1}{N^{|x-y|}(\lambda +h_{|x-y|})}
+(N-1)\sum^\infty_{j=|x-y|+1}\frac{1}{N^j(\lambda +h_j)}.
\end{eqnarray*}
By (\ref{3.2.9}), $\pi_\lambda (x,y)\rightarrow \infty$ as
$\lambda \rightarrow 0$, and by (\ref{3.2.8}),
\begin{eqnarray*}
&&\sum_{y\in \Omega_N}\pi_\lambda (x,y) F(y) = (N-1)F(x)
\sum^\infty_{j=1}{1\over \lambda N^j +NDd_{j-1}}\\
&&+\sum_{y\neq x}F(y)\biggl[-\frac{1}{\lambda N^{|x-y|} +NDd_{|x-y|-1}}
+(N-1)\sum^\infty_{j=|x-y|+1}\frac{1}{\lambda N^j +NDd_{j-1}}\biggr].
\end{eqnarray*}

We know from the proof of Proposition \ref{prop3.3.1} with $\mu=1$ that
$$
\sum_j \frac{1}{\lambda N^j  +\;\hbox{\rm const}\; d_j}\sim\;\hbox{\rm const}\;
\sum_{j=1}^{-\log\lambda/\log N}
\frac{1}{d_j}\;\;\hbox{\rm as}\;\;\lambda\rightarrow 0,
$$
where the right-hand side is slowly varying as $\lambda \rightarrow 0$. The result then follows from Theorem 1 of \cite{[8]}. 
\end{proof}
\medskip
\begin{remark}\label{rem3.4.2}\rm
(a)
In the case of  $d$-dimensional simple symmetric random walks, for $d=1$ the norming is
$t^{1/2}$ and the limit is the truncated normal distribution, and for
$d=2$ the norming is $\log t$ and the limit is the exponential distribution
\cite{[8]}. Hence, form the point of view of occupation time  the critically recurrent
random walks in Proposition \ref{prop3.4.1} behave like $2$-dimensional simple symmetric
random walks.

\noindent
(b) Recall that the recurrent $c^j$-random walk with $c<1$ behaves differently from the  random walks above (Remark \ref{rem3.3.3}(a)). In particular, in contrast with Proposition 3.4.1 the continuous time $c^j$-random walk with $c<1$ does not satisfy an occupation time result as above. Indeed, condition (A) of \cite{[8]} is satisfied with the norming $g(\lambda)=\sum_j1/(\lambda N^j +\;\hbox{\rm const}\;c^j)$ (denoted by $h(s)$ in \cite{[8]}), and by Theorem 2 of \cite{[8]}, if there existed an occupation time limit distribution as $t\rightarrow\infty$, then $g(\lambda)$ would necessarily be of the form $g(\lambda)=\lambda^{-\alpha}L(\lambda^{-1})$ for some $\alpha,0<\alpha\leq 1$, and slowly varying $L(\lambda^{-1})$, and by a Tauberian theorem we would have $G_t\sim t^\alpha L(t)/\Gamma(\alpha+1)$ as $t\rightarrow\infty$. But it is shown in \cite{[9]} (Lemma 3.1.1) that $G_t\sim$ const $t^{-\gamma}h_t$ where $\gamma$ is the degree (\ref{3.2.3}) $(-1<\gamma<0)$, and $h_t$ is the function 
$$h_t=\sum^\infty_{j=-\infty}(b a^{j-1}t)^\gamma(1-
e^{-b a^{j-1}t}),\,\,\,t>0,$$
where  $a=c/N$ and $b=(N^2-c)/N(N-1)$, 
and this function is slowly oscillating but not slowly varying.
\end{remark}

\subsection {Distance Markov chain}\label{dmc}
\setcounter{equation}{0}
Some properties of random walks on $\Omega_N$ depend only on the distance from $0$, which we  study in this subsection. This is more easily done in discrete time.
We exemplify with   the $c^j$-random walk (with $\mu =1$ for simplicity)
to show  explicit results.
\begin{definition}\label{def3.5.1} \rm
Let
$\{\xi_n\}$ be the  $r_j$-random walk on $\Omega_N$ defined by (\ref{3.1.1}) and let
\begin{equation}\label{eq:3.5.1}
Z_n=|\xi_n|.
\end{equation}
$\{Z_n\}$ is a Markov chain on $\ene_0=\{0,1,2,\ldots\}$ called {\it distance Markov chain}.
\end{definition}
 
 We denote the transition probability of $\{Z_n\}$ by
$p_{ij}=\pee [Z_{n+1}=j|Z_n=i]$.

\noindent
\begin{proposition} \label{prop3.5.2} The transition probabilities $p_{ij}$ are as given as follows:
\vglue .25cm
\noindent
{\it  (1) $r_j$-random walk:}
\begin{equation}\label{eq:3.5.2}
\begin{array}{lcll}
p_{ij}&=&r_j, &j>i,\nonumber\\[.10cm]
&& r_1+\cdots +r_{i-1}+r_i\displaystyle{{N-2\over N-1}}=
1-\displaystyle{{r_i\over N-1}}-\sum^\infty_{j=i+1}r_j,\quad &j=i(\neq 0),\quad (p_{00}=0),\nonumber\\[.15cm]
& &r_i\displaystyle{{1\over N^{i-j}}},\quad&0<j<i,\nonumber\\[.15cm]
&&r_i\displaystyle{{1\over N^{i-1}(N-1)}},\quad &0=j<i.
\end{array}
\end{equation}
\noindent
(2) $c^j$-random walk:

\noindent

\begin{equation}\label{eq:3.5.3}
\begin{array}{lcll}
p_{ij}  &=&\left(1-\displaystyle{{c\over N}}\right)
\left(\displaystyle{{c\over N}}\right)^{j-1},&j>i,\nonumber\\[.25cm]
&&1-\left(\displaystyle{{c\over N}}\right)^i-\displaystyle{{1-c/N\over N-1}}
\left(\displaystyle{{c\over N}}\right)^{i-1}
=1-\left(\displaystyle{{c\over N}}\right)^i
\displaystyle{{N-2\over N-1}}-\left(\displaystyle{{ c\over N}}\right)^{i-1}
\displaystyle{{1\over N-1}},& j=i(\neq 0),\nonumber\\[.25cm]
&&\left(1-\displaystyle{{c\over N}}\right)\left(\displaystyle{{c\over N^2}}
\right)^{i-1}N^{j-1}
=\left(1-\displaystyle{{c\over N}}\right)\left(\displaystyle{{c\over N}}
\right)^{i-1}
\displaystyle{{1\over N^{i-j}}},& 0<j<i,\nonumber\\[.25cm]
&&\left(1-\displaystyle{{c\over N}}\right)\left(\displaystyle{{c\over N^2}}
\right)^{i-1}\displaystyle{{1\over N-1}},&0=j<i.
\end{array}
\end{equation}
\end{proposition}
\begin{proof}
The proof relies on  the ultrametric property:
  $|x|<|y|\Rightarrow|y-x|=|y|$, and
$|x|=|y|, x\neq y\Rightarrow |y-x|=|y|$. We prove (\ref{eq:3.5.2}):
\vglue.5cm
$j>i$: A jump of $\{Z_n\}$ from $i$ to $j$ is  the same as from 0 to $j$.
\vglue.5cm
$j=i\;\;(\neq 0)$: This happens in two ways:

\begin{itemize}
\item[(i)] for each $k=1,\ldots,i-1$, $\{\xi_n\}$ jumps to a point with the same $l$-coordinates, $l=k+1,\ldots,i$, and different $k$-coordinate as the previous point, which occurs with probability $r_k$, and all such points are favorable, or
\item[(ii)] $\{\xi_n\}$ jumps to a point with $i$-coordinate different from that of the
previous point and from $0$, which occurs with probability $r_i$, and there are $N^{i-1}(N-2)$ favorable possibilities out of $N^{i-1}(N-1)$.
\end{itemize}
\vglue.2cm
$0<j<i$: $\{\xi_n\}$ jumps a distance $i$ from the previous point, which occurs with probability $r_i$, and there are $N^{j-1}(N-1)$ favorable possibilities out of $N^{i-1}(N-1)$.
\vglue.5cm
$0=j<i$: This is as the previous case with one favorable possibility out of
$N^{i-1}(N-1)$.
\vglue.5cm

(\ref{eq:3.5.3}) is immediate from (\ref{eq:3.5.2}).
\end{proof}
\medskip

We next state without proof some elementary results that follow directly from Proposition \ref{prop3.5.2}.
\begin{proposition}\label{prop3.5.3} Let
$\tau_j=\inf\{n:Z_n\geq j\}$,  $j\geq 1$, and $T_i=$ first exit time of $\{Z_n\}$ from $i$ (starting at $i$). Then
\vglue .25cm
\noindent
 (1) $r_j$-random walk:
\begin{eqnarray}
\label{eq:3.5.4}
\pee_0[\tau_j=n]&=&\left(
\sum^{j-1}_{i=1}r_i\right)^{n-1}\sum^\infty_{i=j}r_i,\quad n=1,2,\ldots ,
\quad \ee_0(\tau_j)={1\over \sum^\infty_{i=j}r_i}.\\
\label{eq:3.5.5}
\pee_i[T_i=n]&=&p^{n-1}_{ii}(1-p_{ii}),\quad n=1,\ldots,
\quad \ee_iT_i={1\over 1-p_{ii}}.
\end{eqnarray}

\noindent
 (2)  $c^j$-random walk:
\begin{equation}
\label{eq:3.5.6}
\pee_0[\tau_j=n]=
\left(
1-\left({c\over N}\right)^{j-1}\right)^{n-1}
\left({c\over N}\right)^{j-1},\, n=1,2,\ldots ,
\quad \ee_0(\tau_j)=\left({N\over c}\right)^{j-1}.
\end{equation}
\begin{equation}\label{eq:3.5.7}
\ee_iT_i
=\left({N\over c}\right)^i
{N-1\over N(1+1/c)-2}.
\end{equation}
\end{proposition}

\begin{remark}\label{rem3.5.4} \rm
For the $c^j$-random walk we have
from (\ref{eq:3.5.3}) that $p_{ii}\approx 1$ for large $i$ or large
$N$,
and for every $i$,
$p_{i,i+1}/p_{i,i-1}=c$ and
$\sum^\infty_{j=i+1}p_{ij}/\sum^{i-1}_{j=0}p_{ij}=c(N-1)/(N-c)= 1\, \,({\rm resp.}\, >1,\, <1)\, 
\,\,\,{\rm iff}\quad c=1\, ({\rm resp.} >1, \, <1)$. It is interesting that these quotients are independent of $i$.
This shows that the  walk tends to stay at the same distance from $0$
and the value of $c$ determines the tendency to go away from  or towards 0.
(\ref{eq:3.5.7}) shows that
 the  walk  stays at distance $i$ an average of the order of
$(N/c)^i$ steps before making a jump to another distance.
Since $\sum^i_{j=0}p_{ij}=((N/c)^i-1)\sum^\infty_{j=i+1}p_{ij}$, in one step from $i$ the distance chain is $(N/c)^i-1$ times more likely to stay within distance $i$ from $0$ than it is to jump to a larger distance from $0$.
%
%\noindent
%(b) For $c\leq 1, \{\xi_n\}$ is null recurrent since it is an irreducible recurrent random walk on an infinite countable Abelian group, hence
%$\{Z_n\}$ is also null recurrent.
%For $c>1,\{Z_n\}$ is transient since  $\{\xi_n\}$ is transient. In the recurrent case $Z_n$ has a unique invariant measure $\pi$ such that $\pi(0)=1$. \marginpar{maybe here we don't need a reference}
%The results on transience/recurrence and $k$-strong/weak transience of $\{Z_n\}$ are the same as those of (the continuous-time version of) $\{\xi_n\}$ by Proposition \ref{prop2.2}, which are given by (\ref{growthg}).  Note that it
% would be  harder to obtain these results  for $\{Z_n\}$ by means of a
%discrete time formula like (\ref{gzetaform})
%because it involves the powers of the complicated transition matrix \ref{eq:3.5.3}.
\end{remark}
\begin{proposition}\label{prop3.5.5} For the $c^j$-random walk, consider the expected distance from $0$ of $\{Z_n\}$ after one step
 starting from $i$, i.e.,  $D_i=\sum\limits^\infty_{j=1}jp_{ij}$. We have
\begin{eqnarray}
D_0&=&{N\over N-c},\nonumber\\ \label{eq:3.5.8}
D_i&=&i+\left({c\over N}\right)^{i-1}\left[{c\over N-c}-{(N-c)(N^i-1)\over N^i(N-1)^2}\right],\qquad i>0
\end{eqnarray}
{\it For}
$c=1$,
\begin{equation}\label{eq:3.5.9}
D_i=i+{1\over N^{2i-1}(N-1)},\quad i>0.
\end{equation}
\end{proposition}
\begin{corollary} \label{cor3.5.6}

\noindent
(1) For $c\geq 1$, $D_i> i$ for all $i$.

\noindent
 (2) For $c<1, D_i<i$  iff
\begin{equation}\label{eq:3.5.10}
i>L_N(c):=
{1\over \log N}\biggl(-\log\biggl(1-c\biggl({N-1\over N-c}\biggr)^2\biggr)
\biggr).
\end{equation}
\end{corollary}
{\bf Proofs} of Proposition \ref{prop3.5.5} and Corollary \ref{cor3.5.6}:
The  calculations use (\ref{eq:3.5.3}) and the standard summation formulas

\begin{eqnarray*}
\sum^n_{j=1}jx^j&=&{x-(n+1)x^{n+1}+nx^{n+2}\over (1-x)^2},\\
\sum^\infty_{j=n}jx^j&=&{nx^n-(n-1)x^{n+1}\over (1-x)^2},\quad 0<x<1.
\end{eqnarray*}
For $i=0$:
\begin{eqnarray*}
D_0&=&
\left(1-{c\over N}\right)\sum^\infty_{j=1}j\left({c\over N}\right)^{j-1}
=
\left(1-{c\over N}\right){N\over c}{c/N\over (1-c/N)^2}=
{1\over 1-c/N}=
{N\over N-c}.
\end{eqnarray*}
For $i>0$:
\begin{eqnarray*}
D_i&=&
i\left[1-\left({c\over N}\right)^i-{N-c\over N(N-1)}\left({c\over N}
\right)^{i-1}\right]\\
&&+\left(1-{c\over N}\right)\sum^\infty_{j=i+1}j\left({c\over N}\right)
^{j-1}+\left(1-{c\over
N}\right)\left({c\over N^2}\right)^{i-1}\sum^{i-1}_{j=1}jN^{j-1}\\
&=&
i\left[1-\left({c\over N}\right)^i-{N-c\over N(N-1)}\left({c\over N}
\right)^{i-1}\right]\\
&+&\left(1-{c\over N}\right){N\over c}{(i+1)(c/N)^{i+1}-i(c/N)^{i+2}\over
(1-c/N)^2}\\
&+&\left(1-{c\over N}\right)\left({c\over N^2}\right)^{i-1}{1\over
N}{N-iN^i+(i-1)N^{i+1}\over (N-1)^2}\\
&=&i+\left({c\over N}\right)^{i-1}
\left[{c\over N-c}-{(N-c)(N^i-1)\over N^i(N-1)^2}\right].\\
\end{eqnarray*}

\noindent
The  term in square brackets is equal to
$$
{cN^i(N-1)^2-(N-c)^2(N^i-1)\over (N-c)N^i(N-1)^2},
$$
and the numerator equals $N^i[c(N-1)^2-(N-c)^2]+(N-c)^2$, which is positive for all $i$ iff $c(N-1)^2\geq (N-c)^2$, iff $c\geq 1$. Hence
for $c\geq 1$, $D_i> i$
for all $i$.

For $c<1, D_i<i$ iff $N^i[(N-c)^2-c(N-1)^2]>(N-c)^2$,
iff
\vglue.3cm
\noindent
$\hfill i>\displaystyle{{1\over \log N}} \log \displaystyle{{(N-c)^2\over
(N-c)^2-c(N-1)^2}}=L_N(c). \hfill\Box$

\begin{remark}\label{rem3.5.7} \rm
(a) Since for $c\geq 1$ (i.e. for non-negative degree of the walk) the drift is positive, in this case $\{Z_n\}$ is a submartingale. For $c<1$,
$ \{Z_n\}$ behaves like a submartingale for $i\leq L_N(c)$,
and when it exceeds $L_N(c)$ it stops behaving that way because the drift becomes negative. Note that $L_N(c)\rightarrow\infty$ as $c\nearrow 1$. In the case of (Euclidean) $d$-dimensional 
Brownian motion (i.e., $c=N^{1-2/d}$, see Remark \ref{rem3.2.2}(a)),
 $\{Z_n\}$ is the analogue of  a Bessel process,
but  Bessel processes do not behave the way $\{Z_n\}$ does. This exhibits a qualitative difference between hierarchical random walks
 and Euclidean processes, which is due to the ultrametric structure of $\Omega_N$.

\noindent
(b) For $c<1$, let
$$T_N(c)=
\lfloor L_N{(c)}
\rfloor+1,$$
\vglue .25cm
\noindent
i.e., $\tau_{T_N{(c)}}$ is  the  time of the first jump over the threshold $L_N(c)$ where the drift of $\{Z_n\}$ becomes negative. Then, from (3.5.6) and (3.5.10),
$$\lim_{N\to \infty} \ee_0(\tau_{T_N(c)})={1\over 1-c}.$$
\end{remark}

We give next some results on the  maximal process $Z^*_n: 
=\max\limits_{1\leq m\leq n} Z_m$, $n=1,2,\ldots $.
\begin{proposition}
\label{prop3.5.8} For $j\geq 1$,

\noindent
(1)  $r_j$-random walk:
\begin{equation}\label{eq:3.5.11}
\pee_0\left[Z_n^* = j\right]=\left(\sum^j_{i=1}r_i\right)^n-\left(\sum^{j-1}_{i=1}r_i\right)^n,
\end{equation}
\begin{equation}\label{eq:3.5.12}
\pee_0\left[Z_n^*\geq j\right]=
1-\left(\sum^{j-1}_{i=1}r_i\right)^n.
\end{equation}
\noindent
(2)  $c^j$-random walk:
\begin{eqnarray}\label{eq:3.5.13}
\pee_0\left[Z_n^*=j\right]
&=&\left(1-\left(
{c\over N}\right)
^j\right)^n-
\left(
1-\left(
{c\over N}\right)^{j-1}\right)^n,\\\label{eq:3.5.14}
\pee_0\left[Z_n^*\geq
j\right]&=&1-\left(1-\left({c\over N}\right)^{j-1}\right)^n.
\end{eqnarray}
\end{proposition}
\begin{proof}

\noindent
(1)
\begin{eqnarray*}
 \pee_0[Z^*_n\geq j]&=&\pee_0[Z^*_n\geq j,Z^*_{n-1}\geq j]+\pee_0[Z^*_n\geq j,
Z^*_{n-1}<j]\\
&=&\pee_0[Z^*_{n-1}\geq j]+\pee_0[Z^*_n\geq j, Z^*_{n-1}<j].
\end{eqnarray*}
By (\ref{eq:3.5.4}),
$$
\pee_0[Z^*_n\geq j,Z^*_{n-1}<j]=
\pee_0[\tau_j=n]
=\left(\sum^{j-1}_{i=1}r_i\right)^{n-1}\sum^\infty_{i=j}r_i,
$$
so
$$
\pee_0[Z^*_n\geq j]=\pee_0[Z^*_{n-1}\geq j]+
\left(\sum^{j-1}_{i=1}r_i\right)^{n-1}
\sum^\infty_{i=j}r_i,$$
hence
$$
\pee_0[Z^*_n\geq j]=\sum^\infty_{i=j}r_i\sum^{n-1}_{\ell=0}
\left(\sum^{j-1}_{i=1}r_i\right)^\ell=1-\left(\sum^{j-1}_{i=1}r_i\right)^n,
$$
and (\ref{eq:3.5.11}) follows.

\noindent
(2) (\ref{eq:3.5.13}) and (\ref{eq:3.5.14}) are special cases of
(\ref{eq:3.5.11}) and (\ref{eq:3.5.12}). 
\end{proof}
\smallskip

The next corollaries are easy consequences (see Remark \ref{rem3.1.5} for Corollary \ref{cor3.5.10}).
\begin{corollary}\label{cor3.5.9}For $j\geq 1$,

\noindent
(1)  $r_j$-random walk:
\begin{equation}\label{eq:3.5.15}
\pee_0\left[Z^*_n=j\right] \sim \left(\sum^j_{i=1}r_i \right)^{n}\quad as \quad n\rightarrow \infty .
\end{equation}
\noindent
(2) $c^j$-random walk:
\begin{equation}\label{eq:3.5.16}
\pee_0\left[Z^*_n=j\right]\sim \left( 1-\left({c\over N}\right)^j\right)^n\quad\mbox{\it as}\quad n\rightarrow \infty.
\end{equation}
\end{corollary}

\begin{corollary}\label{cor3.5.10} For the $(\mu, (\eta^j), N)$-random walk with $\mu \ge 1$,
\vglue.5cm
\noindent
(1) 
\begin{equation}\label{eq:3.5.17}
\lim\limits_{j\rightarrow \infty}\pee_0\left[Z^*_{\lfloor N^{j/\mu}\rfloor}\leq j\right]=\left\{
\begin{array}{c}
0\\
1/e\\
1
\end{array}\right.\;\;\; \hbox{\it iff} \;\;\; \eta \left\{
\begin{array}{ll}
>&  1,\\
= & 1,\\
< & 1
\end{array}\right.
\end{equation}

\vglue.5cm
\noindent
(2)  For $j\geq 1$,
\begin{equation}\label{eq:3.5.18}
\lim_{N\rightarrow\infty}\pee_0
\left[Z^*_{\lfloor N^{j/ \mu}\rfloor}=\ell\right]=
\left\{
\begin{array}{ll}
e^{-\eta^j},&\ell=j,\\
1-e^{-\eta ^j},&\ell=j+1.
\end{array}\right.
\end{equation}
\end{corollary}

\begin{remark}\label{rem3.5.11}\rm
  Corollary \ref{cor3.5.10}  shows  that  $N^{j/\mu}$ is  the right time scale for observing the exit behaviour of a $(\mu, (\eta^j), N)$-random walk from a closed ball of radius $j$. Asymptotically as $N\rightarrow \infty$, only the closed ball of radius $j$ and the surrounding closed ball of radius $j+1$ are relevant. In \cite{[10]} we consider the cases
$\mu = 1,2$ and we study the behaviour of  branching systems on  a  sequence of nested closed balls of increasing radii in $\Omega_N$, which  due to the 
behaviour just described lead to  separation of time scales (see also  Remark \ref{rem3.2.12}(b)) and, as a consequence, to a cascade of quasiequilibria as $N\rightarrow\infty$.
\end{remark}

The following result explains  why it is easier to compute  probabilities for   $Z_n^*$ than for $Z_n$.

\begin{proposition} \label{prop3.5.12} $Z^*_n, n=1,2, \ldots$  is a Markov chain with transition matrix $Q=(q_{ij})$ given by
\begin{equation}\label{eq:3.5.19}
q_{ij}=\left\{
\begin{array}{ll}
0, &j<i,\nonumber\\
\sum^i_{k=1}r_k, &j=i,\\
r_j, &j>i,\nonumber
\end{array}\right.
\end{equation}
\end{proposition}
\begin{proof} Assume  $Z^*_n=i$. Then
$$
Z^*_{n+1} =\left\{
\begin{array}{lcl}
i & \hbox{\rm iff} & |\eta_{n+1}|\leq i,\\
i+k, \;\; k\geq 1 & \hbox{\rm iff} & |\eta_{n+1}|=i+k,
\end{array}\right.
$$
where $\eta_{n+1}$ is the $(n+1)st$ step of the random walk $\{\xi_n\}$,
independently of $Z_1, \ldots, Z_n$. Then the form of $Q$ is obvious.
\end{proof}
\smallskip

\begin{remark}\label{rem3.5.13} \rm
Proposition \ref{prop3.5.12} reflects the fact that  in an ultrametric space all interior points of a closed ball are at the ``center''. Clearly, 
Euclidean random walks do not have the property in this proposition because it matters where inside the ball the  jump starts from. However, it is worthwhile to mention a behaviour of  simple random walk on $\mathbb Z^2$ which has certain features of separation of time scales, with  close connections to the Erd\"os-Taylor theorem (see \cite{CG86} and references therein): Consider the ball $B_R$ with radius $R$ centered around the origin. For all $0 < a < a'$, and large $t$, the walk starting in $x \in B_{t^{a/2}}$ is at time $t^{a'}$ ``nearly uniformly'' distributed on  $B_{t^{a'/2}}$, independently of the starting position.   
\end{remark}

We now give some  results on the moments of $Z_n$ and the 
rate of escape for of the $c^j$-random walk.

\begin{proposition}\label{prop3.5.14} (1) For the $c^j$-random walk and  for all $n\geq 1$ and any $M>0$ ($M$ not necessarily an integer),
\begin{equation}\label{eq:3.5.20} \ee_0 (Z^*_n)^M=\sum^\infty_{j=1}j^M \biggl(\frac{c}{N}\biggr)^j \biggl(\frac{N}{c}-1\biggr)\sum\limits^n_{k=1}\biggl(1-\biggl(\frac{c}{N}\biggr)^j\biggr)^{n-k}\biggl(1-\biggl(\frac{c}{N}\biggr)^{j-1}\biggr)^{k-1},
\end{equation}
and
\begin{equation}\label{eq:3.5.21}
\ee_0 Z^M_n \leq n\frac{N-c}{c}\sum^\infty_{j=1}j^M\biggl(\frac{c}{N}\biggr)^j \biggl(1-\biggl(\frac{c}{N}\biggr)^j\biggr)^{n-1},
\end{equation}
(2)
\begin{equation}\label{eq:3.5.22}
\lim_{n\rightarrow \infty} \frac{1}{n}\ee_0 Z^M_n=0.
\end{equation}
(3)
\begin{equation}
\label{eq:3.5.23}
\lim_{n\rightarrow \infty} \frac{Z_n}{n}=0\quad \hbox{\it a.s.}
\end{equation}
\end{proposition}
\begin{proof}(1) Let $a=c/N$. By (3.5.13),
\begin{eqnarray*}
\ee_0 (Z^*_n)^M&=& \sum^\infty_{j=1}j^M\left[(1-a^j)^n -(1-a^{j-1})^n \right]\\
&=& \sum^\infty_{j=1}j^M(a^{j-1}-a^j)\sum^n_{k=1}(1-a^j)^{n-k}(1-a^{j-1})^{k-1},
\end{eqnarray*}
which is (\ref{eq:3.5.20}).

To obtain  (\ref{eq:3.5.21}) we use the obvious inequalities
$$
\ee_0 Z^M_n \leq \ee_0 (Z^*_n)^M,
$$
and
\begin{eqnarray*}
(a^{j-1}-a^j) \sum^n_{k=1} (1-a^j)^{n-k} (1-a^{j-1})^{k-1}&=&(a^{-1}-1)a^j (1-a)^{n-1}\sum^n_{k=1}\left(\frac{1-a^{j-1}}{1-a^j}\right)^{k-1}\\
&\leq  &n(a^{-1}-1)a^j (1-a^j )^{n-1} .
\end{eqnarray*}
(2) (\ref{eq:3.5.22}) follows from (\ref{eq:3.5.21}) by dominated convergence.
\vglue .25cm
\noindent
(3) (\ref{eq:3.5.23}) follows from (\ref{eq:3.5.21}) by Chebyshev's inequality and the Borel-Cantelli lemma. 
\end{proof}
\smallskip
\begin{remark}\label{rem3.5.15} \rm
(a) The  transition matrix (\ref{eq:3.5.19}) of $Z^*_n$ for the $c^j$-random walk is
$$
q_{ij}=\left\{
\begin{array}{ll}
0, & j<i\\
1-(c/N)^i, & j=i,\\
(1-c/N) (c/N)^{j-1}, & j>i,
\end{array}\right.
$$
and the $n$-step transition matrix $Q^n=(q^{(n)}_{ij})$ is given by
\begin{equation}\label{eq:3.5.24}
q^{(n)}_{ij}=\left\{
\begin{array}{ll}
0, & j<i,\\
(1-(c/N)^i)^n, & j=i,\\
(1-(c/N)^j)^n-(1-(c/N)^{j-1})^n, & j>i.
\end{array}\right.
\end{equation}
\noindent
(b) Proposition \ref{prop3.5.14} (3) means that the  rate of escape of the $c^j$-random walk is $0$. The following result, which is more precise than (\ref{eq:3.5.23}), is obtained using (\ref{eq:3.5.24}):
\begin{equation}
\mathbb P_0[Z_n^\ast \ge \delta \log n] \sim  \frac{\lfloor\delta \log(c/N)\rfloor}{1+\lfloor\delta \log(c/N)\rfloor}n^{1+\lfloor\delta \log(c/N)\rfloor} \quad \mbox{ as  } n \to \infty 
\end{equation}
for all $\delta > 1/\log(N/c)$, and this implies for any $\delta > 2/\log(N/c)$,
\begin{equation}
\mathbb P_0[Z_n^\ast \ge \delta \log n \,\,\,{\rm  i.\,o.}] = 0. 
\end{equation}
\end{remark}

\section{Occupation time fluctuations of  branching systems}\label{sec4}
\setcounter{equation}{0}

In this section we apply the results on  the operator $G_t$  obtained in subsection \ref{sec3.3}  to derive asymptotic results for the occupation time fluctuations of branching systems. To keep the presentation self-contained, we first give a short review of the subject.

\subsection{Incomplete potentials and growth functions}\label{sec4.1} 

Multilevel branching systems were introduced by Dawson and Hochberg \cite{[15]} and they have been studied by several authors \cite{[9], [10], [16], [18], [24], [26], [27], [57]}. In  addition to the individual particle branching there is an independent branching of families of related particles (2-level branching), and this idea  can be extended to higher levels of branching. The main difficulty in dealing with these models is that the independence of behaviour of individual particles no longer holds due to the higher-level branchings.

Here we assume that the group $S$ is locally compact with countable base, Haar measure $\rho$, and the process $X$ has stationary independent increments which are symmetric and have a strictly positive density with respect to $\rho$. In the analysis of large time  occupation time fluctuations of $k$-level
 branching particle systems on  $S$ (where $k=0$ corresponds to absence of branching), a basic problem  consists in finding a norming $a_t$ such that the occupation time fluctuation
\begin{equation}
\label{eq:3.1}
{1\over a_t}\int^t_0 ({\cal X}_s-\ee {\cal X}_s)ds
\end{equation}
has a non-trivial limit in distribution as $t\rightarrow \infty$, where ${\cal X}_s$ in the empirical measure of the particle system at time $s$.
Under appropiate assumptions on the system (suitable initial conditions, critical binary branchings), it turns out that
$\ee {\cal X}_t=\rho$  for all $t$, and  in the cases of recurrent and of $k$-weakly transient motion the form of $a_t$ is dictated by the order of the growth of operator $G_t$ defined by (\ref{Gt}) and its powers  as $t\rightarrow \infty$. Precisely, $a_t$ is determined by $G_t$ for recurrent motion, by $G^2_t$ for weakly transient motion, and by $G^3_t$ (or $G^2_tG$) for $2$-weakly transient motion.

Occupation time fluctuation limits of up to $2$-level branching systems were investigated in \cite{[9]}, to which we refer the reader for more information and details. For the $0$-level and the $1$-level particle systems the initial condition was taken to be a Poisson random field with intensity $\rho$. The $1$-level system has a ``Poisson-type'' equilibrium state, and for the $2$-level system the initial condition was taken to be a Poisson random field of ``$2$-level particles'' whose intensity is the canonical measure of the equilibrium state of the $1$-level system. The moments of this canonical measure involve the potential operator $G$ \cite{[9]} (Appendix), and this implies that one has to deal with  $G^2_tG$ rather than $G^3_t$ (e.g. (\ref{3.3.12})). A different initial condition that can be assumed for the $2$-level system is a Poisson random field with intensity measure $\delta_{\delta_x}\rho (dx)$, and  this would lead to  $G^3_t$ in place of
$G^2_tG$. (In case $R_t$ defined in Remark \ref{rem2.10} (a) decreases like $t^{-\gamma}$ for some $\gamma > 0$, then $G_t^2G$ and $G_t^3$ have the same order of growth, see \cite{[9]}, Lemma 2.4.2.)

%but it makes the proof of Proposition 3.3.1 for $\mu=3$ more elaborate.

It is shown in \cite{[9]} that for each $k$-level branching system, if the growth of $G_t, G^2_t$, etc., is given by an increasing function $f_t$, then the norming $a_t$ for the occupation time fluctuation (\ref{eq:3.1}) is
\begin{equation}
\label{eq:3.3}
a_t =\biggl(\int^t_0 f_s ds\biggr)^{1/2}.
\end{equation}
For $k$-strongly transient motion $a_t$ is the ``classical'' noming $a_t=t^{1/2}$.

For the $\alpha$-stable process on $\erre^d$ (with no branching),
$$
a_t=\left\{
\begin{array}{ll}
t^{1-d/2\alpha}&{\rm for}\quad \alpha>d,\\
(t\log t)^{1/2}&{\rm for} \quad \alpha=d,\\
t^{1/2}&{\rm for}\quad \alpha <d.
\end{array}\right.
$$
Note that $t^{1-d/2\alpha}\rightarrow t^{1/2}$ as $\alpha \searrow d$, so there is a discontinuity in the order of the growth at $\alpha=d$, and for this value of $\alpha$ the  ``critical'' fluctuations of the occupation time are bigger than $t^{1/2}$. The critical case corresponds to $\gamma=0$, where $\gamma$ is the degree of the $\alpha$-stable process given by (\ref{eq:2.19}).

For Brownian motion $(\alpha =2)$ on $\erre^d$
and the $0$-level system (no branching):
$$
a_t=\left\{
\begin{array}{lclcl}
t^{3/4} & \hbox{\rm for}&  d=1,\\
(t\log t)^{1/2} &  \hbox{\rm for} & d=2,\\
t^{1/2} & \hbox{\rm for} & d\geq 3,
\end{array}\right.
$$
\cite{[5], [17]}. The same pattern is repeated for the $1$-level branching system (individual particle branching) $2$ dimensions higher \cite{CG85}, where the critical case corresponds to $\gamma = 1$, and for the $2$-level branching system (individual branching and  family branching) $4$ dimensions higher \cite{[9]}, where the critical case corresponds to $\gamma = 2$.

In the general setting of branching systems on locally compact Abelian groups the $t\rightarrow \infty$ limits of the occupation time fluctuations are Gaussian random fields described in detail in \cite{[9]}. The Gaussian property is due to the finiteness of the variance of the branching laws. A class of infinite variances branching laws leads to stable random fields \cite{[9]}.

\subsection{Occupation time fluctuations of $j^\beta$-branching random walks}\label{sec3.6}

The occupation time fluctuation limits of branching systems of $c^j$-random walks are given  in \cite{[9]}.
A different situation occurs for the class of hierarchical random walks in subsection \ref{sec3.3}.    For illustration  we consider the
 $j^\beta$-random walk (Example \ref{ex3.2.6}, $d_j =(j+1)^\beta$, $\beta \geq 0$).
We obtain the following result from (\ref{3.3.12}) and (\ref{eq:3.3}) for
$\mu=1$ ($0$-level system), $\mu=2$ ($1$-level system) and
$\mu=3$ ($2$-level system):

$$
a_t = \left\{
\begin{array}{ll}
t^{1/2}(\log t)^{(1-\mu \beta )/2}&{\rm for}\quad \beta<1/\mu,\\
(t\log\log t)^{1/2}&{\rm for}\quad \beta=1/\mu,\\
t^{1/2}&{\rm for}\quad \beta>1/\mu.
\end{array}\right.
$$
The  forms of the limit Gaussian random fields of the occupation time fluctuations can be  obtained from \cite{[9]} (Theorems 2.2.1 to 2.2.3), and the  constants can be computed from (\ref{3.3.12}). For example, for the $0$-level system with transient motion $(\beta >1, a_t =t^{1/2})$ the covariance kernel of the limit Gaussian  field, obtained from (\ref{3.1.5}), is
$$
k(x,y) =\frac{2N}{D} \biggl[(N-1) \zeta (\beta )+ \biggl(\delta_{0, |x-y|}-1\biggr) |x-y|^{-\beta}-(N-1)\sum^{|x-y|}_{j=1}j^{-\beta}\biggr],
$$
where $D$ is the normalizing constant in (3.1.8) and $\zeta (\cdot)$ is the Riemann Zeta function.
For the 1-level system  in the critical case $(\beta=1/2,a_t=(t\log\log t)^{1/2})$, the covariance kernel of the limit Gaussian field is a constant 
$(=(N-1)/ND^2)$, hence the occupation time fluctuation limits in all regions of $\Omega_N$ are perfectly correlated.

\vglue .5cm
\noindent
{\bf Acknowledgments}
\vglue .5cm
The authors thank the referee for the careful reading of the paper, and comments and suggestions that led to significant improvements. They also thank the hospitality of The Fields Institute (Toronto, Canada), Carleton University (Ottawa, Canada), the Center for Mathematical Research (CIMAT, Guanajuato, Mexico),  the Johann Wolfgang Goethe University (Frankfurt, Germany), and the Erwin Schr\"odinger Institute (Vienna, Austria), where mutual working visits took place. L.G.G. also thanks the Institute of Mathematics, National University of Mexico (UNAM), where he spent a sabbatical during 2002.

\end{document}